\pdfoutput=1
\documentclass[10pt,review,authoryear,times]{elsarticle}

\makeatletter
\def\ps@pprintTitle{%
   \let\@oddhead\@empty
   \let\@evenhead\@empty
   \let\@oddfoot\@empty
   \let\@evenfoot\@oddfoot
}
\makeatother

\usepackage{setspace}
\usepackage{xcolor}
\definecolor{elsblue}{RGB}{0, 128, 172}
\usepackage[colorlinks=true]{hyperref}
\usepackage{hyperref}
\usepackage{xparse}
\usepackage{amsmath}
\usepackage{amsfonts}

\newcommand{\bxi}{\boldsymbol{\xi}}

\newcommand{\dVblank}{\, {\rm d}V}
\newcommand{\dvblank}{\, {\rm d}v}
\newcommand{\dxi}{\, {\rm d}\bxi}

\newcommand{\eps}{\varepsilon}

\NewDocumentCommand \dV{ o }{%
    \IfNoValueTF{#1}{\dVblank}%
    {
        \dV_{#1}
    }
}

\NewDocumentCommand \dv{ o }{%
    \IfNoValueTF{#1}{\dvblank}%
    {
        \dv_{#1}
    }
}

\NewDocumentCommand \mc{ m }{%
    \mathcal{#1}%
}

\NewDocumentCommand \p{ m m }{%
    \frac{\partial #1}{\partial #2}%
}

\NewDocumentCommand \Dp{ m m }{%
    \dfrac{\partial #1}{\partial #2}%
}

\NewDocumentCommand \an{ m }{%
    \langle {#1} \rangle%
}

\NewDocumentCommand \mbf{ m }{%
    \mathbf{#1}
}

\NewDocumentCommand \mbb{ m }{%
    \mathbb{#1}
}

\NewDocumentCommand \bbar{ m }{%
    \bar{\mbf{#1}}
}

\NewDocumentCommand \bs{ m }{%
    \boldsymbol{#1}
}

\NewDocumentCommand \stateOne{ m }{
    \underline{\mbf{#1}}%
}

\NewDocumentCommand \stateTwo{ m o }{%
    \IfNoValueTF{#2}{\stateOne{#1}}
    {
    \stateOne{#1}\an{#2}%
    }
}

\NewDocumentCommand \s{ m o o }{%
    \IfNoValueTF{#3}{\stateTwo{#1}[#2]}%
    {
    \underline{\mbf{#1}}(#3)\an{#2}%
    }
}

\NewDocumentCommand \stateOneI{ m m }{
    \underline{#1}_{#2}%
}

\NewDocumentCommand \stateTwoI{ m m o }{%
    \IfNoValueTF{#3}{\stateOneI{#1}{#2}}
    {
    \stateOneI{#1}{#2}\an{#3}%
    }
}

\NewDocumentCommand \sI{ m m o o }{%
    \IfNoValueTF{#4}{\stateTwoI{#1}{#2}[#3]}%
    {
        \stateOneI{#1}{#2}(\mbf{#4})\an{#3}
    }
}

\NewDocumentCommand \stateOneN{ m }{
    \underline{#1}%
}

\NewDocumentCommand \stateTwoN{ m o }{%
    \IfNoValueTF{#2}{\stateOneN{#1}}
    {
    \stateOneN{#1}\an{#2}%
    }
}

\NewDocumentCommand \sN{ m o o }{%
    \IfNoValueTF{#3}{\stateTwoN{#1}[#2]}%
    {
    \underline{#1}(#3)\an{#2}%
    }
}

\NewDocumentCommand \scs{ m }{
    \underline{#1}%
}

\NewDocumentCommand\presuper{ m m }{%
  {\mathop{}%
   \mathopen{\vphantom{#2}}^{#1}%
   \kern-1\scriptspace%
   #2}
}

\NewDocumentCommand\presuperz{ m m m }{%
  {\mathop{}%
   \mathopen{\vphantom{#2}}^{#1}%
   \kern-3\scriptspace%
   #2\kern-8\scriptspace_{#3}%
}
}

\usepackage{pgfplots}

\usepackage[english]{babel}
\usepackage[utf8x]{inputenc}
\usepackage[T1]{fontenc}

\usepackage[a4paper,top=0.65in,bottom=0.65in,left=0.6in,right=0.6in,marginparwidth=0.5in]{geometry}

\usepackage{amsmath}
\usepackage{graphicx}
\usepackage{subfig}

\usepackage{stfloats}

\usepackage[nameinlink]{cleveref}
\crefname{equation}{Eq.}{Eqs.}
\crefname{figure}{Fig.}{Figs.}
\crefname{section}{Section}{Sections}
\crefname{remark}{Remark}{Remarks}
\crefname{appendix}{}{}

\newtheorem{proposition}{Proposition}
\newdefinition{remark}{Remark}

\makeatletter
\newcommand*\bigcdot{\mathpalette\bigcdot@{.8}}
\newcommand*\bigcdot@[2]{\mathbin{\vcenter{\hbox{\scalebox{#2}{$\m@th#1\bullet$}}}}}
\makeatother

\newcommand{\I}{{\rm I}}
\newcommand{\J}{{\rm J}}

\makeatletter
\newcommand{\vast}{\bBigg@{3}}
\newcommand{\Vast}{\bBigg@{3.5}}
\newcommand{\vasT}{\bBigg@{4}}
\makeatother

\begin{document}

\hypersetup{allcolors=elsblue}

\setlength{\abovedisplayskip}{6pt}
\setlength{\belowdisplayskip}{6pt}

\begin{frontmatter}

\title{On the stability of the generalized, finite deformation correspondence model of peridynamics}


\author[mymainaddress]{Masoud Behzadinasab\corref{mycorrespondingauthor}}
\ead{behzadi@utexas.edu}

\author[mymainaddress,mysecondaryaddress]{John T. Foster}
\cortext[mycorrespondingauthor]{Corresponding author.}


\address[mymainaddress]{Department of Aerospace Engineering \& Engineering Mechanics, The University of Texas at Austin, United States}
\address[mysecondaryaddress]{Hildebrand Department of Petroleum \& Geosystems Engineering, The University of Texas at Austin, United States}

\begin{abstract}
    A class of peridynamic material models known as \emph{constitutive correspondence models} provide a bridge between classical continuum mechanics and peridynamics.  These models are useful because they allow well-established local constitutive theories to be used within the nonlocal framework of peridynamics. A recent finite deformation correspondence theory \citep{foster2018generalized} was developed and reported to improve stability properties of the original correspondence model \citep{silling2007peridynamic}.  This paper presents a stability analysis that indicates the reported advantages of the new theory were overestimated. Homogeneous deformations are analyzed and shown to exibit unstable material behavior at the continuum level. Additionally, the effects of a particle discretization on the stability of the model are reported. Numerical examples demonstrate the large errors induced by the unstable behavior.  Stabilization strategies and practical applications of the new finite deformation model are discussed.
\end{abstract}


\end{frontmatter}


\section{Introduction}
\label{sec:intro}
Peridynamics was introduced by \citet{silling2000reformulation} in order to represent the mechanics of continuous and discontinuous media using a single consistent set of equations. The theory has gained notable attention, due to its capabilities in naturally representing material failure, removing the need for complicated numerical treatments in dealing with discontinuities, i.e.\ cracks. Several peridynamic material models have been proposed since introduction of the theory \citep{silling2007peridynamic,foster2010viscoplasticity,mitchell2011nonlocal,foster2018generalized} and applied to numerous problems, from micro to macro scales, mostly involving brittle fracture \citep{gerstle2007peridynamic,askari2008peridynamics,ha2010studies,bobaru2018intraply,jafarzadeh2018peridynamic,behzadinasab2018perturbation}. A special class of material models was proposed by \citet{silling2007peridynamic}, called \emph{constitutive correspondence models}, which provide a convenient framework for using classical constitutive relations in peridynamic material modeling. 

While the correspondence model approach offers convenience, it has shortcomings associated with material instability. Stability issues have been associated with the continuum model and its discretization \citep{littlewood2011nonlocal,tupek2014extended,silling2017stability}. Zero-energy mode oscillations appear in meshless simulations that use correspondence material models and result in severe degradation in the accuracy of displacement fields, particularly in areas of steep gradients \citep{silling2017stability}. These oscillations where originally reported as numerical issues; however, \citet{tupek2014extended} demonstrated some kinematically unphysical deformation modes that are permitted in the mathematical formulation of the model. To address the issue, they introduced a nonlocal family of Seth-Hill nonlinear bond strain measures to create what they called an ``extended constitutive correspondence formulation''. Recently, \citet{foster2018generalized} generalized Tupek's model and presented a finite deformation correspondence model which does not suffer from any {\em surface effects} -- a well-known issue in peridynamics \citep{mitchell2015position}. 

Presented here is an analysis of the stability, in the sense of Lyapunov, of the generalized finite deformation model and its meshless discretization. It is shown that the new constitutive correspondence theory exhibits zero-energy mode oscillations as a manifestation of a new type of material instability which occurs even in the presence of uniform deformations. The root cause of the problem is investigated to better understand the role of nonlocality on the Seth-Hill strain measures. 

The rest of the paper is outlined as follows: In \cref{sec:background} a brief background of the peridynamic theory and an overview of the correspondence material modeling is presented. In \cref{sec:homogeneous}  the stability of equilibrium for homogeneous deformations is analyzed at the continuum level. \cref{sec:discrete} presents the effects of particle discretization on the stability. A set of unstable equilibrium points are analytically found in \cref{sec:homogeneous,sec:discrete}, limited to pure hydrostatic deformations. Numerical examples are provided in \cref{sec:numerical} to demonstrate the stability of the theory in a more general loading scenario. In \cref{sec:discussion} the connection between the stability method presented in this work and the method of \citet{silling2017stability} is discussed along with possible stabilization strategies and the source of instability. \cref{sec:conclusions} offers concluding remarks.

\section{Background on peridynamics and correspondence modeling}
\label{sec:background}
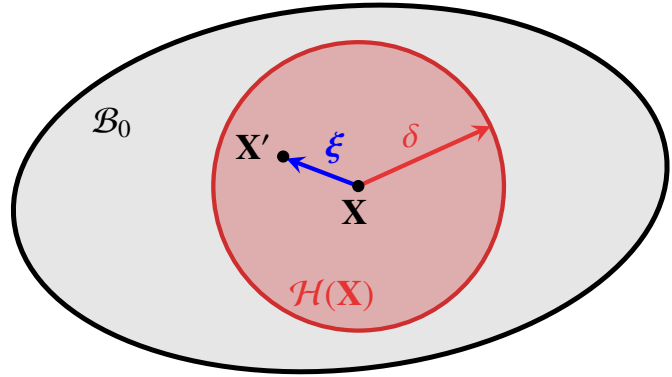
\begin{figure}[!htbp]
  \centering
  \usetikzlibrary{arrows}
\definecolor{wwqqzz}{rgb}{0.9,0.2,0.2}
\resizebox {0.5\columnwidth} {!} {
\begin{tikzpicture}[line cap=round,line join=round,>=triangle 45,x=1.0cm,y=1.0cm]
\clip(2.425951811078324,1.585) rectangle (8.922494932630698,5.286381878454878);
\draw [line width=1.2pt,color=wwqqzz,fill=wwqqzz,fill opacity=0.3] (5.849948764351535,3.44010247129693) circle (1.4212208160764792cm);
\draw [rotate around={5.60353836762431:(5.67,3.42)},line width=1.4pt,fill=black,fill opacity=0.1] (5.67,3.42) ellipse (3.2072177542608573cm and 1.7877767543085596cm);
\draw [color=wwqqzz](5.0524255998835725,2.6499593252980906) node[anchor=north west] {$\mathcal{H}(\mathbf{X}) $};
\draw (3.1223653156857765,4.341249265059049) node[anchor=north west] {$\mathcal{B}_0$};
\draw [->,>=stealth,line width=1.2pt,color=wwqqzz] (5.846585714244864,3.4433732823330048) -- (7.141666113705931,4.032837570161464);
\draw (5.559761346163395,3.42550182388164) node[anchor=north west] {$\mathbf{\mathbf{X}}$};
\draw [color=wwqqzz](6.154310290131715,4.185402972004443) node[anchor=north west] {$\mathbf{\delta}$};
\draw (4.520167611931382,4.082171506731418) node[anchor=north west] {$\mathbf{\mathbf{X}'}$};
\draw [color=blue](5.40786059132693,4.102017799786023) node[anchor=north west] {$\bs{\xi}$};
\draw [->,>=stealth,line width=1.2pt,color=blue] (5.846585714244864,3.4433732823330048) -- (5.141964479047388,3.7144771109123407);
\begin{scriptsize}
\draw [fill=black] (5.846585714244864,3.4433732823330048) circle (1.5pt);
\draw [fill=black] (5.111964479047388,3.7344771109123407) circle (1.5pt);
\end{scriptsize}
\end{tikzpicture}
}
  \caption{Schematic of nonlocality within a peridynamic body with horizon $\delta$.}
  \label{fig:PD}
\end{figure}

The peridynamic theory recasts the mechanics of solid deforamtions in a nonlocal setting. The mechanics of continuous and discontinuous bodies are treated using the same set of equations, which fundamentally allows for discontinuities in the deformation field. Suppose $\mbf{x}(\mbf{X}, t)$ is the current position of a material point $\mbf{X}$ in a deformed peridynamic body $\mc{B}$. The motion of $\mbf x$ is influenced by other material points through long-range interactions with a \emph{nonlocal family}. The family of $\mbf{X}$ is defined in its undeformed configuration by
\begin{equation*}
  \mc{H}(\mbf{X}) = \left\{ \mbf{X} + \bxi \ \vert \ \mbf{X} + \bxi \in \mc{B}_0, \, 0 < \vert\bxi\vert \leq \delta \right\} ,
\end{equation*}
where the vector $\bxi$ is called a {\em bond}, and $\delta$ is the radius of a spherical neighborhood, called the {\em horizon}. Note that the point $\mbf{X}$ itself is excluded from the family. \cref{fig:PD} shows an example of a peridynamic body.

Peridynamics describes the balance of linear momentum using the following integro-differential equation:
\begin{equation}
  \rho_0\ddot{\mbf{u}}(\mbf{X}, t) = \int_{\mc{H}(\mbf{X})} \left[ \s{T}[\bxi][\mbf{X}, t] - \s{T}[-\bxi][\mbf{X}+\bxi, t] \right] \dxi + \rho_0 \, \mbf{b}(\mbf{X}, t) ,
  \label{eqn:PDEOM}
\end{equation}
where $\rho_0$ is the reference mass density, $\mbf{u}$ is the displacement, and $\mbf{b}$ is a prescribed body force density field. $\s{T}$ is a vector field (with units of force per square volume), called the \textit{force vector state}, which is used to describe the constitutive interactions of material points. 

\citet{silling2007peridynamic} introduced the idea of \emph{constitutive correspondence} where a nonlocal analog of a kinematic variable is utilized as a bridge between the classical (i.e.\ local) and peridynamic theories. \citet{silling2007peridynamic} defined a nonlocal deformation gradient tensor by
\begin{equation}
  \bbar{F} = \left( \int_{\mc{H}} \s{\omega}[\bxi] \, \s{Y}[\bxi] \otimes \bxi \dxi \right) \mbf{K}^{-1} ,
  \label{eqn:F}
\end{equation}
where $\s{\omega}[\bxi]$ is called the {\em influence state}, describing the relative degree of interactions between the neighboring points; $\s{Y}[\bxi]$ is the deformed image of the bond in its current configuration,
\begin{equation*}
  \s{Y}[\bxi][\mbf{X},t] = \mbf{x}(\mbf{X}+\bxi,t) - \mbf{x}(\mbf{X},t);
\end{equation*}
and $\mbf{K}$ is a symmetric second-order shape tensor, given by
\begin{equation*}
  \mbf{K} = \int_{\mc{H}} \s{\omega}[\bxi] \, \bxi \otimes \bxi \dxi .
  \label{eqn:K}
\end{equation*}
The authors showed that in a continuous media, the peridynamic deformation gradient recovers \emph{exactly} the local quantity under homogeneous deformations. Using this kinematic variable as a metric and incorporating classical theories by means of the first Piola-Kirchhoff stress tensor $\bs{\sigma} = \bs{\sigma}(\bbar{F})$, the force vector state was derived
\begin{equation}
  \sN{T}_k = \s{\omega}[\bs{\xi}] \, \sigma_{kp} K^{-1}_{pq} \xi_q .
  \label{eqn:TF}
\end{equation}
It was observed \citep{tupek2014extended,silling2017stability} that the nonlocal analog $\bbar{F}$ allows for unphysical deformation modes such as matter interpenetration and can result in loss of material stability. The instability can effectively destroy the accuracy of displacement fields in areas of large displacement gradients. The peridynamic finite deformation correspondence framework \citep{tupek2014extended} and its generalization \citep{foster2018generalized} were proposed to improve the stability of the original correspondence theory by using a family of nonlinear strain measures. As an intermediate step in determining the force vector-state, the generalized theory \citep{foster2018generalized} utilizes a family of nonlocal Lagrangian Seth–Hill strains. A set of nonlocal right Cauchy-Green deformation tensors is introduced 
\begin{equation}
  \bbar{C}_{(m)} = \left[\int_{\mc{H}} \s{\omega}[\bxi] \, \Bigg(\frac{\left|\s{Y}[\bxi]\right|}{\left|\bxi\right|}\Bigg)^{2m} \frac{\bxi \otimes \bxi}{\left|\bxi\right|^2}\dxi\right]:\mbb{L}^{-1} , \quad m\neq0 ,
  \label{eqn:C}
\end{equation}
where $\mbb{L}$ is a symmetric fourth-order shape tensor defined by
\begin{equation*}
  \mbb{L} = \int_{\mc{H}} \s{\omega}[\bxi] \frac{\bxi \otimes \bxi \otimes \bxi \otimes \bxi}{\left|\bxi\right|^4}\dxi .
  \label{eqn:L}
\end{equation*}
The peridynamic Lagrangian Seth–Hill strains are then
\begin{equation}
  \bbar{E}_{(m)} = \frac{1}{2m}\left(\bbar{C}_{m}-\mbf{I}\right) , \quad m\neq0 ,
  \label{eqn:strains}
\end{equation}
and 
\begin{equation*}
  \bbar{E}_{(0)} = \left[\int_{\mc{H}} \s{\omega}[\bxi] \ln\Bigg(\frac{\left|\s{Y}[\bxi]\right|}{\left|\bxi\right|}\Bigg) \frac{\bxi \otimes \bxi}{\left|\bxi\right|^2}\dxi\right]:\mbb{L}^{-1} ,
\end{equation*}
which is shown to recover exactly their local analogs under uniform deformations. Note that $m=0$ and $m=1$ are two important members of the class, which correspond to the widely used Hencky and Green-Lagrange strain measures. A local hyperelastic material model can then be integrated into peridynamics, where the force vector-state is given as
\begin{equation}
  \sN{T}_k = \s{\omega}[\bs{\xi}] \, S_{(m)ij} \frac{\xi_p \, \xi_q}{|\bs{\xi}|^2} L^{-1}_{pqij} \left(\frac{|\s{Y}|}{|\bs{\xi}|}\right)^{2m} \frac{Y_k}{|\s{Y}|^2},
  \label{eqn:TC}
\end{equation}
where $\mbf{S}_{(m)}$ is the generalized Kirchhoff stress and is a function of $\bbar{E}_{(m)}$. 

\section{Stability under homogeneous deformation}
\label{sec:homogeneous}
In this section  the stability of the generalized, finite deformation correspondence theory when dealing with uniform deformations is analyzed in a continuum setting. It is first shown that homogeneous deformations constitute a set of equilibria for a peridynamic material under static conditions and in the absence of any external body forces. The stability of these equilibrium points in the context of the finite deformation model are discussed.

\subsection{Static equilibrium}
\label{sec:staticEquilibrium}
For a static problem in the absence of body forces, the equation of motion, \cref{eqn:PDEOM}, reduces to
\begin{equation}
    \int_{\mc{H}(\mbf{X})} \left[ \s{T}[\bxi][\mbf{X}] - \s{T}[-\bxi][\mbf{X} + \bxi] \right] \dxi = 0 .
  \label{eqn:PDstatic}
\end{equation}

For a static homogeneous deformation field with $\mbf{F}(\mbf{X}) = \mbf{F}$, the deformation vector state for a peridynamic bond can be described as
\begin{equation*}
  \s{Y} = \mbf{F} \cdot \bs{\xi} .
\end{equation*}
Then \cref{eqn:TC} simplifies to
\begin{equation}
  \sN{T}_k = \s{\omega}[\bs{\xi}] \, S_{(m)ij} L^{-1}_{pqij} F_{kl} \left(\frac{|\mbf{F}\cdot\bs{\xi}|}{|\bs{\xi}|}\right)^{2m-2} \frac{\xi_p \, \xi_q \, \xi_l}{|\bs{\xi}|^4}.
  \label{eqn:uniformT}
\end{equation}
Limiting the influence state to be spherical, i.e. $\s{\omega}\an{\bs{\xi}} = \s{\omega}_s(|\bs{\xi}|)$, $\s{T}[\bs{\bxi}]$ is clearly an odd function of $\bs{\xi}$. Therefore, for two neighboring points in the bulk of material and away from the boundaries (having the same shape tensor)
\begin{equation}
  \s{T}[\bs{\xi}][\mbf{X}] = -\s{T}[-\bs{\xi}][\mbf{X}+\bxi] .
  \label{eqn:TforTwo}
\end{equation}
Using \cref{eqn:PDstatic,eqn:TforTwo} 
\begin{equation*}
\int_{\mc{H}(\mbf{X})} \left[ \s{T}[\bxi][\mbf{X}] - \s{T}[\bxi][\mbf{X} + \bxi] \right] \dxi = 2 \int_{\mc{H}(\mbf{X})} \s{T}[\bxi][\mbf{X}] \dxi.
\end{equation*}
Since $\s{T}[\bs{\xi}]$ is an odd function of $\bs{\xi}$, this integral is zero for a full neighborhood and therefore satisfies \cref{eqn:PDstatic} automatically in the bulk of material. With appropriate boundary conditions, which equilibrate the internal forces in the vicinity of the boundaries, static equilibrium can be achieved for any homogeneous deformation. A static uniform deformation cannot be obtained if any body force is present in the bulk of material.

\subsection{Stability}
\label{sec:stability}
Analyzing the Jacobian matrix of a mechanical system is a way to understand its stability at equilibrium. If all the eigenvalues of the Jacobian matrix are negative, the system is stable.  Presence of one or more positive eigenvalues indicates unstable behavior. A zero eigenvalue can correspond to either stable or unstable behavior. 

Let the net force exerted on $\mbf{X}$ in \cref{eqn:PDEOM} be called $\mbf{f}(\mbf{X})$. In statics, it is given by
\begin{equation*}
  \mbf{f}(\mbf{X}) = \int_{\mc{H}(\mbf{X})} \left[ \s{T}[\bxi][\mbf{X}] - \s{T}[-\bxi][\mbf{X} + \bxi] \right] \dxi .
\end{equation*}
Using partial differentiation, the Jacobian matrix is computed as 
\begin{align*}
  \sN{J}_{kl}(\mbf X)\an{\bs\eta} = \p{f_k(\mbf{X})}{x_l(\mbf{X} + \bs\eta)} ,
  \label{eqn:J}
\end{align*}
or
\begin{equation}
  \p{f_k(\mbf X)}{x_l(\mbf X + \bs\eta)} = \int_{\mc{H}(\mbf X)} \left( \p{\scs{T}_{k}(\mbf X)\an\bxi}{x_l(\mbf{X+\bs\eta})} - \p{\scs{T}_{k}(\mbf X + \bxi)\an{-\bxi}}{x_l(\mbf X + \bs\eta)} \right) \dxi .
  \label{eqn:df}
\end{equation}
Rewriting \cref{eqn:TC}
\begin{equation*}
    \scs{T}_k(\mbf X)\an\bxi = \scs\omega\an\bxi \frac{\xi_p \, \xi_q}{\vert\bs{\xi}\vert^{2m+2}} L^{-1}_{pqij} \, S_{(m)ij}(\mbf{X}) \, \vert\scs{Y}_n\an\bxi \, \scs{Y}_n\an\bxi\vert^{m-1} \scs{Y}_k\an\bxi ,
\end{equation*}
and taking the partial derivative 
\begin{align}
\p{\sN{T}_{k}(\mbf X)\an\bxi}{x_l(\mbf X + \bs\eta)} 
  = \scs\omega\an\bxi \frac{\xi_p \, \xi_q}{\vert\bs{\xi}\vert^{2m+2}} L^{-1}_{pqij} \, \Bigg[ &\p{S_{(m)ij}(\mbf{X})}{x_l(\mbf X + \bs\eta)} \vert\s{Y}[\bxi]\vert^{2m-2} \sN{Y}_k\an\bxi 
   + S_{(m)ij}(\mbf{X}) \, (m-1)\bigg(2\sN{Y}_n\an{\bxi}\p{\sN{Y}_n\an{\bxi}}{x_l(\mbf X + \bs\eta)}\bigg)\left|\s{Y}[\bxi]\right|^{2m-4} \sN{Y}_k\an\bxi \nonumber \\
   &+ S_{(m)ij}(\mbf{X}) \left|\s{Y}[\bxi]\right|^{2m-2} \p{\sN{Y}_k\an\bxi}{x_l(\mbf X + \bs\eta)} \Bigg] .
  \label{eqn:dT}
\end{align}
The partial derivative of the deformation state is calculated 
\begin{align}
    \p{\scs{Y}_k\an\bxi}{x_l(\mbf X + \bs\eta)} &= \p{\left(x_k(\mbf X + \bxi) - x_k(\mbf X) \right)}{x_l(\mbf X + \bs\eta)} , \nonumber \\
                                                &= \underline{\delta}\langle\bxi,\bs\eta\rangle \, \delta_{kl} - \underline{\delta}\langle\mbf 0,\bs\eta\rangle \, \delta_{kl} ,
  \label{eqn:dY}
\end{align}
where $\underline{\delta}\langle\bxi,\bs\eta\rangle$ is the Kronecker-delta state, i.e.
\begin{equation*}
    \underline{\delta}\langle\bxi,\bs\eta\rangle = 
    \begin{cases}
        1 & \mbox{if} \quad \bxi = \bs\eta, \\
        0 & \mbox{otherwise}.
    \end{cases}
\end{equation*}
Using $\mbf{S}_{(m)}=\mbf{S}(\mbf{\bar{E}}_{(m)})$ and the chain rule
\begin{equation}
  \p{S_{(m)ij}(\mbf X)}{x_l(\mbf X + \bs\eta)} = \p{S_{(m)ij}}{\bar{E}_{(m)rs}} \p{\bar{E}_{(m)rs}(\mbf X)}{x_l(\mbf X + \bs\eta)}.
  \label{eqn:dS}
\end{equation}
Taking the derivative of \cref{eqn:C},
\begin{align}
    &\p{\bar{C}_{(m)rs}(\mbf X)}{x_l(\mbf X + \bs\eta)} 
    = L^{-1}_{turs} \int_{\mc{H}(\mbf X)} \scs\omega\an\bxi \, m \left(2\scs{Y}_n\an\bxi\p{\scs{Y}_n\an\bxi}{x_l(\mbf X + \bs\eta)}\right) \vert\s{Y}[\bxi]\vert^{2m-2} \frac{\xi_t \, \xi_u}{\vert\bxi\vert^{2m+2}} \dxi , \nonumber \\
    &~~~~= L^{-1}_{turs} \int_{\mc{H}(\mbf X)} 2m \, \scs\omega\an\bxi \, \bigg(\sN{Y}_n\an\bxi \left(\scs\delta\an{\bxi, \bs\eta} - \scs\delta\an{\mbf{0}, \bs\eta}\right)\delta_{nl} \bigg) \, \vert\s{Y}[\bxi]\vert^{2m-2} \frac{\xi_t \, \xi_u}{\vert\bxi\vert^{2m+2}} \dxi , \nonumber \\
    &~~~~= 2m \, L^{-1}_{turs} \Bigg( \scs\omega\an{\bs\eta} \, \sN{Y}_l\an{\bs\eta} \, \vert\s{Y}\an{\bs\eta}\vert^{2m-2} \frac{\eta_t \, \eta_u}{\vert\bs\eta\vert^{2m+2}} \, {\rm d}\bs\eta 
    - \scs\delta\an{\mbf 0, \bs \eta} \int_{\mc{H}(\mbf X)} \scs\omega\an\bxi \, \sN{Y}_l \, \an\bxi \, \vert\s{Y}[\bxi]\vert^{2m-2} \frac{\xi_t \, \xi_u}{\vert\bxi\vert^{2m+2}} \dxi \Bigg).
    \label{eqn:partialC}
\end{align}
Using \cref{eqn:strains,eqn:partialC},
\begin{align*}
    \p{\bar{E}_{(m)rs}(\mbf X)}{x_l(\mbf X + \bs\eta)} &= \frac{1}{2m} \p{\bar{C}_{(m)rs}(\mbf X)}{x_l(\mbf X + \bs\eta)} \nonumber , \\
    &= L^{-1}_{turs} \Bigg( \scs\omega\an{\bs\eta} \, \sN{Y}_l\an{\bs\eta} \, \vert\s{Y}\an{\bs\eta}\vert^{2m-2} \frac{\eta_t \, \eta_u}{\vert\bs\eta\vert^{2m+2}} \, {\rm d}\bs\eta - \scs\delta\an{\mbf 0, \bs \eta} \int_{\mc{H}(\mbf X)} \scs\omega\an\bxi \, \sN{Y}_l\an\bxi \, \vert\s{Y}[\bxi]\vert^{2m-2} \frac{\xi_t \, \xi_u}{\vert\bxi\vert^{2m+2}} \dxi \Bigg).
\end{align*}
It can be shown that this equation is valid for $m=0$. Similar to the argument made in deriving \cref{eqn:uniformT}, for a uniform deformation, the integrand is odd, and the integral is equal to zero. Therefore,
\begin{equation}
  \p{\bar{E}_{(m)rs}(\mbf X)}{x_l(\mbf X + \bs\eta)} = L^{-1}_{turs} \scs\omega\an{\bs\eta} \, \sN{Y}_l\an{\bs\eta} \, \vert\s{Y}\an{\bs\eta}\vert^{2m-2} \frac{\eta_t \, \eta_u}{\vert\bs\eta\vert^{2m+2}} \, {\rm d}\bs\eta .
  \label{eqn:dE}
\end{equation}
Substitute \cref{eqn:dY,eqn:dS,eqn:dE} into \cref{eqn:dT} and simplify to get 
\begin{align}
  \p{\sN{T}_{k}(\mbf X)\an\bxi}{x_l(\mbf X + \bs\eta)} = \
  \s{\omega}[\bxi] \frac{\xi_p \, \xi_q}{\left|\bs{\xi}\right|^{2m+2}} L^{-1}_{pqij} \left|\s{Y}[\bxi]\right|^{2m-2} \Bigg[ & \p{S_{(m)ij}}{\bar{E}_{(m)rs}} L^{-1}_{turs} \, \sN{Y}_k\an{\bxi} \, \s{\omega}[\bs\eta] \, \sN{Y}_l\an{\bs\eta} \left|\s{Y}[\bs\eta]\right|^{2m-2} \frac{\eta_t \, \eta_u}{\left|\bs\eta\right|^{2m+2}} \, {\rm d}\bs\eta\nonumber \\ 
  & + S_{(m)ij}(\mbf{X}) \Big( \underline{\delta}\langle\bxi,\bs\eta\rangle - \underline{\delta}\langle\mbf{0},\bs\eta\rangle \Big) \bigg( 2(m-1)\frac{\sN{Y}_k\an{\bxi} \, \sN{Y}_l\an{\bxi}}{\left|\s{Y}[\bxi]\right|^2} + \delta_{kl} \bigg) \Bigg] .
 \label{eqn:TCI}
\end{align}
Noting that $\s{Y}[\bxi] = -\s{Y}[-\bxi]$, $\mbf{S}_{(m)}(\mbf{X}) = \mbf{S}_{(m)}(\mbf{X}+\bxi) = \mbf{S}_{(m)}$, and $\s{\omega}[\bxi] = \s{\omega}[-\bxi]$, we use \cref{eqn:TCI} to derive
\begin{align}
  \p{\sN{T}_{k}(\mbf X + \bxi)\an{-\bxi}}{x_l(\mbf X + \bs\eta)} =
  - \, \s{\omega}[\bxi] \frac{\xi_p \, \xi_q}{\left|\bs{\xi}\right|^{2m+2}} L^{-1}_{pqij} \left|\s{Y}[\bxi]\right|^{2m-2} \Bigg[ & \p{S_{(m)ij}}{\bar{E}_{(m)rs}} L^{-1}_{turs} \, \sN{Y}_k\an{\bxi} \, \s{\omega}[\bs\eta-\bxi] \, \sN{Y}_l\an{\bs\eta-\bxi} \, \left|\s{Y}[\bs\eta-\bxi]\right|^{2m-2} \frac{(\eta_t-\xi_t) (\eta_u-\xi_u)}{\left|\bs\eta-\bxi\right|^{2m+2}} \, {\rm d}\bs\eta \nonumber \\ 
  & + S_{(m)ij} \left( \underline{\delta}\langle\bxi,\bs\eta\rangle - \underline{\delta}\langle\mbf{0},\bs\eta\rangle \right) \bigg( 2(m-1)\frac{\sN{Y}_k\an{\bxi} \, \sN{Y}_l\an{\bxi}}{\left|\s{Y}[\bxi]\right|^2} + \delta_{kl} \bigg) \Bigg] .
 \label{eqn:TCJ}
\end{align}
Substitute \cref{eqn:TCI,eqn:TCJ} in \cref{eqn:df} and rearrange terms to get
\begingroup
\allowdisplaybreaks
\begin{align*}
  \p{f_k(\mbf X)}{x_l(\mbf X + \bs\eta)} = & \, \mkern-10mu \int_{\mc{H}(\mbf{X})} \mkern-10mu \s{\omega}[\bxi] \frac{\xi_p \, \xi_q}{\left|\bs{\xi}\right|^{2m+2}} L^{-1}_{pqij} \left|\s{Y}[\bxi]\right|^{2m-2} \Bigg[ \p{S_{(m)ij}}{\bar{E}_{(m)rs}} L^{-1}_{turs} \, \sN{Y}_k\an{\bxi} \, \s{\omega}[\bs\eta] \, \sN{Y}_l\an{\bs\eta} \, \left|\s{Y}[\bs\eta]\right|^{2m-2} \frac{\eta_t \, \eta_u}{\left|\bs\eta\right|^{2m+2}} \, {\rm d} \bs \eta \Bigg] \, {\rm d} \bxi \notag \\
  & + \mkern-10mu \int_{\mc{H}(\mbf{X})} \mkern-10mu \s{\omega}[\bxi] \frac{\xi_p \, \xi_q}{\left|\bs{\xi}\right|^{2m+2}} L^{-1}_{pqij} \left|\s{Y}[\bxi]\right|^{2m-2} \Bigg[ \p{S_{(m)ij}}{\bar{E}_{(m)rs}} L^{-1}_{turs} \, \sN{Y}_k\an{\bxi} \, \s{\omega}[\bs\eta-\bxi] \, \sN{Y}_l\an{\bs\eta-\bxi} \left|\s{Y}[\bs\eta-\bxi]\right|^{2m-2} \frac{(\eta_t-\xi_t) (\eta_u-\xi_u)}{\left|\bs\eta-\bxi\right|^{2m+2}} \, {\rm d} \bs \eta \Bigg] \, {\rm d} \bxi \notag \\
  & + \mkern-10mu \int_{\mc{H}(\mbf{X})} \mkern-10mu \s{\omega}[\bxi] \frac{\xi_p \, \xi_q}{\left|\bs{\xi}\right|^{2m+2}} L^{-1}_{pqij} \left|\s{Y}[\bxi]\right|^{2m-2} \notag \Bigg[ 2 S_{(m)ij} \left(\underline{\delta}\langle\bxi,\bs\eta\rangle - \underline{\delta}\langle\mbf{0},\bs\eta\rangle\right) \bigg( 2(m-1) \frac{\sN{Y}_k\an{\bxi} \, \sN{Y}_l\an{\bxi}}{\left|\s{Y}[\bxi]\right|^2} + \delta_{kl} \bigg) \Bigg] \, {\rm d} \bxi 
\end{align*}
\endgroup%
Bringing the constants out of the integrations and grouping terms to further simplify 
\begin{align*}
  \p{f_k(\mbf X)}{x_l(\mbf X + \bs\eta)} = & \p{S_{(m)ij}}{\bar{E}_{(m)rs}} L^{-1}_{turs} L^{-1}_{pqij} \vast\{ \s{\omega}[\bs\eta] \frac{\left|\s{Y}[\bs\eta]\right|^{2m-2} \sN{Y}_l\an{\bs\eta} \, \eta_t \, \eta_u}{\left|\bs\eta\right|^{2m+2}} \, {\rm d} \bs \eta \int_{\mc{H}(\mbf{X})} \s{\omega}[\bxi] \frac{\left|\s{Y}[\bxi]\right|^{2m-2} \sN{Y}_k\an{\bxi} \, \xi_p \, \xi_q}{\left|\bs{\xi}\right|^{2m+2}} \dxi \nonumber \\ 
  & \qquad \qquad \qquad \quad + \, {\rm d} \bs \eta \int_{\mc{H}(\mbf{X})} \Bigg[ \s{\omega}[\bxi] \, \s{\omega}[\bs\eta-\bxi] \frac{\left|\s{Y}[\bxi]\right|^{2m-2} \sN{Y}_k\an{\bxi} \, \xi_p \, \xi_q}{\left|\bs{\xi}\right|^{2m+2}} \frac{\left|\s{Y}[\bs\eta-\bxi]\right|^{2m-2} \sN{Y}_l\an{\bs\eta-\bxi} \, (\eta_t-\xi_t) (\eta_u-\xi_u)}{\left|\bs\eta-\bxi\right|^{2m+2}} \Bigg] \dxi \vast\} \nonumber \\ 
  & + 2 S_{(m)ij} L^{-1}_{pqij} \, \s{\omega}[\bs\eta] \frac{\left|\s{Y}[\bs\eta]\right|^{2m-2} \eta_p \, \eta_q}{\left|\bs{\eta}\right|^{2m+2}} \bigg( 2(m-1)\frac{\sN{Y}_k\an{\bs\eta} \, \sN{Y}_l\an{\bs\eta}}{\left|\s{Y}[\bs\eta]\right|^2} + \delta_{kl} \bigg) \, {\rm d} \bs \eta \nonumber \\
  & - 2 S_{(m)ij} L^{-1}_{pqij} \, \underline{\delta}\langle\mbf{0},\bs\eta\rangle \mkern-10mu \int_{\mc{H}(\mbf{X})} \mkern-10mu \s{\omega}[\bxi] \frac{\left|\s{Y}[\bxi]\right|^{2m-2} \xi_p \, \xi_q}{\left|\bs{\xi}\right|^{2m+2}} \bigg[ 2(m-1)\frac{\sN{Y}_k\an{\bxi} \, \sN{Y}_l\an{\bxi}}{\left|\s{Y}[\bxi]\right|^2} + \delta_{kl} \bigg] \dxi . 
\end{align*}%
%
The first integral is zero for uniform deformations. Using \cref{eqn:strains}, the above equation is further reduced%
\begin{align*}
  \p{f_k(\mbf X)}{x_l(\mbf X + \bs\eta)} = & 
  \p{S_{(m)ij}}{\bar{E}_{(m)rs}} L^{-1}_{turs} L^{-1}_{pqij} \, {\rm d} \bs \eta \mkern-8mu \int_{\mc{H}(\mbf{X})} \Bigg[ \s{\omega}[\bxi] \, \s{\omega}[\bs\eta-\bxi] \frac{\left|\s{Y}[\bxi]\right|^{2m-2} \sN{Y}_k\an{\bxi} \, \xi_p \, \xi_q}{\left|\bs{\xi}\right|^{2m+2}} \frac{\left|\s{Y}[\bs\eta-\bxi]\right|^{2m-2} \sN{Y}_l\an{\bs\eta-\bxi} \, (\eta_t-\xi_t) (\eta_u-\xi_u)}{\left|\bs\eta-\bxi\right|^{2m+2}} \Bigg] \, {\rm d} \bxi \nonumber \\
  & + 2 S_{(m)ij} L^{-1}_{pqij} \, \s{\omega}[\bs\eta] \frac{\left|\s{Y}[\bs\eta]\right|^{2m-2}\eta_p \, \eta_q}{\left|\bs{\eta}\right|^{2m+2}} \bigg[ 2(m-1)\frac{\sN{Y}_k\an{\bs\eta} \, \sN{Y}_l\an{\bs\eta}}{\left|\s{Y}[\bs\eta]\right|^2} + \delta_{kl} \bigg] \, {\rm d} \bs \eta \nonumber \\
  & - 2 S_{(m)ij} L^{-1}_{pqij} \, \underline{\delta}\langle\mbf{0},\bs\eta\rangle \int_{\mc{H}(\mbf{X})} \mkern-10mu \s{\omega}[\bxi] \frac{\left|\s{Y}[\bxi]\right|^{2m-2} \xi_p \, \xi_q}{\left|\bs{\xi}\right|^{2m+2}} \bigg[ 2(m-1)\frac{\sN{Y}_k\an{\bxi} \, \sN{Y}_l\an{\bxi}}{\left|\s{Y}[\bxi]\right|^2} + \delta_{kl} \bigg] \, {\rm d} \bxi .
\end{align*}
The first term is nonzero if $\mbf{X}+\bs\eta$ is among the neighbors of neighbors, the second term is nonzero if $\mbf{X}+\bs\eta$ is a neighbor of $\mbf{X}$, and the third term is nonzero only when $\bs\eta=\mbf{0}$. The special case $\bs\eta=\mbf{0}$ and $k=l$ corresponds to the diagonal values of the Jacobian matrix, which need to be negative as a necessary condition for stability. Focusing on the stability of $\mbf{X}$ when itself is perturbed. 
\begin{align}
  \p{f_k(\mbf X)}{x_l(\mbf X)}\Bigg|_{l=k} = &
  - \p{S_{(m)ij}}{\bar{E}_{(m)rs}} L^{-1}_{turs} L^{-1}_{pqij} \, {\rm d} \mbf{X} \int_{\mc{H}(\mbf{X})} \Bigg[ \s{\omega}^2\an\bxi \Bigg(\frac{\left|\s{Y}[\bxi]\right|}{\left|\bs{\xi}\right|}\Bigg)^{4m-2}\frac{\sN{Y}^2_k\an{\bxi}}{\left|\s{Y}[\bxi]\right|^2} \frac{\xi_p \, \xi_q \, \xi_t \, \xi_u}{\left|\bs{\xi}\right|^{6}} \Bigg] \, {\rm d} \bxi \nonumber \\
  & - 2 S_{(m)ij} L^{-1}_{pqij} \int_{\mc{H}(\mbf{X})} \Bigg[ \s{\omega}[\bxi] \Bigg(\frac{\left|\s{Y}[\bxi]\right|}{\left|\bs{\xi}\right|}\Bigg)^{2m-2} \frac{\xi_p \, \xi_q}{\left|\bs{\xi}\right|^{4}} \bigg( 2(m-1)\frac{\sN{Y}^2_k\an{\bxi}}{\left|\s{Y}[\bxi]\right|^2} + 1 \bigg) \Bigg] \, {\rm d} \bxi . 
  \label{eqn:finaldf}
\end{align}
In the continuum limit the second term dominates as the first one is of $O({\rm d} \mbf{X})$. Hence, in a continuous media 
\begin{align}
  \p{f_k(\mbf X)}{x_l(\mbf X)}\Bigg|_{l=k} \approx - 2 S_{(m)ij} L^{-1}_{pqij} \int_{\mc{H}(\mbf{X})} & \s{\omega}[\bxi] \Bigg(\frac{\left|\s{Y}[\bxi]\right|}{\left|\bs{\xi}\right|}\Bigg)^{2m-2} \frac{\xi_p \, \xi_q}{\left|\bs{\xi}\right|^{4}} \bigg[ 2(m-1)\frac{\sN{Y}^2_k\an{\bxi}}{\left|\s{Y}[\bxi]\right|^2} + 1 \bigg] \, {\rm d} \bxi . 
  \label{eqn:dfcontinuum}
\end{align}
This property depends on $m$, from the finite deformation class, and the loading condition. 

\begin{remark}
\label{rmk:0}
Using $\bs\eta = \mbf{0}$ in \cref{eqn:dE}, and noting that $\s{\omega}[\mbf{0}] = 0$,
\begin{equation*}
  \p{\bar{E}_{(m)ij}(\mbf X)}{x_k(\mbf X)} = \mbf{0} , 
\end{equation*}
therefore,
\begin{equation*}
  \p{S_{(m)ij}(\mbf X)}{x_k(\mbf X)} = \mbf{0}.
\end{equation*}
\citet[Section 2.3]{foster2018generalized} claimed that the finite deformation correspondence theory eliminates the center-point instability of the original correspondence model and hence improves it.  However, it is observed that this type of instability is still present under homogeneous deformations. 
\end{remark}

\subsubsection{Hydrostatic deformations}
\label{sec:hydrostatic}
Consider the pure hydrostatic loading
\begin{equation}
  \sN{Y}_k\an{\bxi} = (1+a) \xi_k 
  \label{eqn:hydro}
\end{equation}
with $a>-1$. Using \cref{eqn:hydro,eqn:dfcontinuum}
\begin{equation}
  \p{f_k(\mbf X)}{x_l(\mbf X)}\Bigg|_{l=k} = - 2 S_{(m)ij} L^{-1}_{pqij} \mkern-8mu \int_{\mc{H}(\mbf{X})} \mkern-8mu \s{\omega}[\bxi] \left| 1+a \right|^{2m-2} \frac{\xi_p \, \xi_q}{\left|\bxi\right|^{4}} \bigg[ 2(m-1)\frac{\xi^2_k}{\left|\bxi\right|^2} + 1 \bigg] \, {\rm d} \bxi . 
  \label{eqn:delF}
\end{equation}
Using the following hydrostatic constitutive relation:
\begin{equation}
  S_{(m)ij} = \kappa a \delta_{ij} ,
  \label{eqn:constitutive}
\end{equation}
where $\kappa$ is a positive material constant (a multiple of bulk modulus, depending on the problem dimension). Using \cref{eqn:delF,eqn:constitutive},
\begin{equation}
  \p{f_k(\mbf X)}{x_l(\mbf X)}\Bigg|_{l=k} = - 2 \kappa \left|1+a\right|^{2m-2} a \Gamma(m) , 
  \label{eqn:Gamma}
\end{equation}
where
\begin{equation}
  \Gamma(m) = L^{-1}_{pqii} \int_{\mc{H}(\mbf{X})} \s{\omega}[\bxi]  \frac{\xi_p \, \xi_q}{\left|\bxi\right|^{4}} \bigg[ 2(m-1)\frac{\xi_k^2}{\left|\bxi\right|^2} + 1 \bigg] \, {\rm d} \bxi . 
  \label{eqn:gammaDef}
\end{equation}
A necessary condition for stability is $\Dp{f_k(\mbf X)}{x_k(\mbf X)} < 0 $; thus, 
\begin{equation}
  a \Gamma(m) > 0  
  \label{eqn:gammaCondition}
\end{equation}
is required. 

\begin{proposition}
Depending on the problem dimension, a critical $m_{\rm cr}$ exists such that $\Gamma(m_{\rm cr}) = 0$.
\end{proposition}
\noindent \textbf{Proof.} \cref{eqn:gammaDef} is simplified to find $m_{\rm cr}$ for 1, 2, and 3 dimensions.

\noindent {\em 1D}: $L^{-1} = \dfrac{1}{V_{\s{\omega}}}$, where $V_{\s{\omega}}$ is the weighted volume
\begin{equation*}
  V_{\s{\omega}} = \int_{\mc{H}} \s{\omega}\an{\bxi} \dxi.
\end{equation*}
Then,
\begin{equation*}
  \Gamma(m) = \frac{1}{V_{\s{\omega}}} \int_{\mc{H}(\mbf{x})} \s{\omega}[\bxi]  \frac{1}{\left|\bxi\right|^{2}} \bigg[ 2(m-1) + 1 \bigg] \, {\rm d} \bxi .
\end{equation*}
\begin{equation*}
  \Gamma(m_{\rm cr}) = 0 \Rightarrow m_{\rm cr} = \frac{1}{2}.
\end{equation*}
\noindent {\em 2D}: $L_{pqii}^{-1} = \dfrac{2}{V_{\s{\omega}}} \delta_{pq}$ (see \cref{sec:AppendixA}), 
\begin{equation*}
  \Gamma(m) = \frac{2}{V_{\s{\omega}}} \int_{\mc{H}(\mbf{x})} \s{\omega}[\bxi]  \frac{1}{\left|\bxi\right|^{2}} \bigg[ 2(m-1) \frac{\xi_k^2}{\left|\bxi\right|^2} + 1 \bigg] \, {\rm d} \bxi .
\end{equation*}
$\Gamma(m)$ is then evaluated to be 
\begin{align*}
  \Gamma(m) &= \frac{2}{V_{\s{\omega}}} \int\displaylimits_{0}^{\delta} \int\displaylimits_{0}^{2\pi} \omega(r) \frac{1}{r^2} \bigg[ 2(m-1) \cos^2(\theta) + 1 \bigg] r \, {\rm d} \theta \, {\rm d} r , \\
  &= \frac{4\pi m}{V_{\s{\omega}}} \int\displaylimits_{0}^{\delta} \omega(r) \frac{1}{r} \, {\rm d} r .
\end{align*}

\begin{equation*}
  \Gamma(m) = 0 \Rightarrow m_{\rm cr} = 0.
\end{equation*}
\noindent {\em 3D}: $L_{pqii}^{-1} = \dfrac{3}{V_{\s{\omega}}} \delta_{pq}$ (cf. \citep[Section 3.2]{foster2018generalized}), then performing calculations similar to 2D results in
\begin{equation*}
  \Gamma(m) = \frac{8\pi m}{V_{\s{\omega}}} \int\displaylimits_{0}^{\delta} \omega(r) {\rm d} r .
\end{equation*}
\begin{equation*}
  \Gamma(m_{\rm cr}) = 0 \Rightarrow m_{\rm cr} = 0.
\end{equation*}

\begin{remark}
\label{rmk:1}
It is evident that 
\begin{equation*}
  \begin{cases}
    m > m_{\rm cr} \Rightarrow \Gamma(m) > 0, \\
    m = m_{\rm cr} \Rightarrow \Gamma(m) = 0, \\
    m < m_{\rm cr} \Rightarrow \Gamma(m) < 0. 
  \end{cases}
  \label{eqn:mGamma}
\end{equation*}
\end{remark}
Therefore, using \cref{eqn:gammaCondition,eqn:mGamma}, there are 3 possibilities: 
\begin{itemize}
  \item $m > m_{\rm cr}$: $a<0$ is not favored. Thus, the system is unstable under pure compression loadings. 
  \item $m = m_{\rm cr}$: $a\Gamma(m)=0$ for any loading, and the unstable loadings cannot be determined yet. 
  \item $m < m_{\rm cr}$: If $a>0$, the system shows unstable behavior. In other words, pure expansions are unstable equilibria for the system. 
\end{itemize}
Note that the violation of \cref{eqn:gammaCondition} is indicative of instability. However, even if \cref{eqn:gammaCondition} holds, stability is not guaranteed; the full Jacobian matrix should be considered for proving stability.

\begin{remark}
\label{rmk:2}
It is realized that in the continuum, for $m\neq m_{\rm cr}$, either the slightest compression or tension will cause instability. It will be shown in \cref{sec:discrete} that discretization would help the stability and reduce the unstable regions.
\end{remark}

\subsection{Stability of Silling's correspondence model}
\label{sec:silling}
The same procedure is repeated for the original correspondence theory \citep{silling2007peridynamic}. Similar to \cref{sec:staticEquilibrium}, it can be shown that homogeneous deformations constitute a set of equilibria for the static case, in the absence of body forces. To analyze the stability of such equilibria, the Jacobian matrix elements are computed. Taking derivative of \cref{eqn:TF},
\begin{equation}
  \p{\sN{T}_{k}(\mbf X)\an\bxi}{x_l(\mbf X + \bs\eta)} = \s{\omega}[\bxi] \p{\sigma_{kp}(\mbf X)}{x_l(\mbf X + \bs\eta)} K^{-1}_{pq} \xi_q .
  \label{eqn:TFI}
\end{equation}
Using $\bs\sigma = \bs\sigma(\bbar{F})$,
\begin{equation}
  \p{\sigma_{kp}(\mbf X)}{x_l(\mbf X + \bs\eta)} = \p{\sigma_{kp}}{F_{mn}} \p{\bar{F}_{mn}(\mbf{X})}{x_l(\mbf X + \bs\eta)} . 
  \label{eqn:sigmaI}
\end{equation}
Using \cref{eqn:F,eqn:dY},
\begin{align*}
  \p{\bar{F}_{mn}(\mbf{X})}{x_l(\mbf X + \bs\eta)}
  &= \left[ \int_{\mc{H}(\mbf{X})} \s{\omega}[\bxi] \left(\underline{\delta}\langle\bxi,\bs\eta\rangle - \underline{\delta}\langle\mbf{0},\bs\eta\rangle\right)\delta_{ml} \, \xi_r \, {\rm d}\bxi \right] K^{-1}_{rn} , \nonumber \\
  &= \delta_{ml} \, \s{\omega}[\bs\eta] \, \eta_r \, {\rm d} \bs \eta \, K^{-1}_{rn} - \underline{\delta}\langle\mbf{0},\bs\eta\rangle \, \delta_{ml} \left( \int_{\mc{H}(\mbf{X})} \s{\omega}[\bxi] \, \xi_r {\rm d}\bxi \right) K^{-1}_{rn} .
\end{align*}
The last integral is zero since the integrand is odd; therefore, 
\begin{equation}
  \p{\bar{F}_{mn}(\mbf{X})}{x_l(\mbf X + \bs\eta)} = \delta_{ml} K^{-1}_{rn} \s{\omega}[\bs\eta] \, \eta_r \, {\rm d} \bs \eta .
  \label{eqn:FI}
\end{equation}
Combining \cref{eqn:df,eqn:TFI,eqn:sigmaI,eqn:FI} to obtain the Jacobian element 
\begin{align*}
  \p{f_k(\mbf X)}{x_l(\mbf X + \bs\eta)} = \p{\sigma_{kp}}{F_{ln}} K^{-1}_{pq} K^{-1}_{rn} \, {\rm d} \bs \eta \int_{\mc{H}(\mbf{X})} \s{\omega}[\bxi] \left[ \s{\omega}[\bs\eta] \, \eta_r + \s{\omega}[\bs\eta-\bxi] \, (\eta_r - \xi_r) \right] \xi_q \, {\rm d} \bxi . 
\end{align*}
Again, paying special attention to the case $\bs\eta=\mbf{0}$
\begin{align} 
  \p{f_k(\mbf X)}{x_l(\mbf X)}
  &= \p{\sigma_{kp}}{F_{ln}} K^{-1}_{pq} K^{-1}_{rn} \, {\rm d} \mbf{X} \int_{\mc{H}(\mbf{X})} \s{\omega}[\bxi] \, \s{\omega}[-\bxi] \, (-\xi_r) \, \xi_q \, {\rm d} \bxi , \nonumber \\
  &= - \p{\sigma_{kp}}{F_{ln}} K^{-1}_{pq} K^{-1}_{rn} \, {\rm d} \mbf{X} \int_{\mc{H}(\mbf{X})} \s{\omega}^2\an\bxi \, \xi_r \, \xi_q \, {\rm d} \bxi .
\end{align}
For a typical influence state
\begin{equation}
  \s{\omega}[\bxi] = 
  \begin{cases}
    1 \qquad \text{if } \left|\bxi\right| \leq \delta , \\
    0 \qquad \text{otherwise} ,
  \end{cases}
  \label{eqn:omega}
\end{equation}
the integral resultant equals the shape tensor, which is 
\begin{equation*}
  K_{ij} = \frac{4\pi}{15}\delta^5 \delta_{ij} 
\end{equation*}
in 3D. Hence,
\begin{equation}
  \p{f_k(\mbf X)}{x_l(\mbf X)} = - \frac{15}{4\pi\delta^5} \p{\sigma_{kp}}{F_{lp}} \, {\rm d} \mbf{X} .
  \label{eqn:dfsimplified}
\end{equation}
If the property of positive-definiteness in $\Dp{\sigma_{ij}}{F_{kl}}$ is inherited in the underlying classical model, all the eigenvalues of $\bar{J}_{kl}(\mbf X)\an{\mbf0}$ would be negative and the model will be stable for any homogeneous deformation. 

The main takeaway is that the unstable behavior identified in the generalized, finite deformation model is not present in Silling's model. There are, however, other types of instabilities associated with Silling's model, such as the zero-energy mode oscillations \citep{breitenfeld2014non, silling2017stability, chowdhury2019modified}, which are associated to non-uniform deformations and particle discretization.

\section{Effect of discretization on the stability}
\label{sec:discrete}
The effect of particle dicretization on the stability of the model is studied in this section. In the discrete form the first term in \cref{eqn:finaldf} is no longer negligible. A general uniform loading is considered first, then it is specialized to the hydrostatic loading condition. Since $m_{\rm cr}$ is different in 1D and higher dimensions (cf. \cref{sec:hydrostatic}), 1D and 2D are discussed separately. In the subsequent sections, to facilitate the notation, $\mbf{X}$, $\mbf{X}+\bxi$, and $\bxi$ are referred to as I, J, and IJ indices, respectively.

\subsection{1D}
\label{sec:1D}
Consider a discretized form of \cref{eqn:finaldf} in a one-dimensional setup and simplify:
\begin{align}
  \p{f^\I}{x^\I} = 
  & - \p{S_{(m)}}{\bar{E}_{(m)}} L^{-1} L^{-1}\Bigg[ \sum_{\mbf{X^J}\in\mc{H^I}} \left(\s{\omega}^{\I\J}\right)^2 \left(\frac{|\sN{Y}^{IJ}|}{|\xi^{\I\J}|}\right)^{4m-2} \left(\frac{1}{|\xi^{\I\J}|}\right)^{2} \Delta V^{\J} \Bigg] \Delta V^\I 
   - 2 \left( 2m-1 \right) S_{(m)} \, L^{-1} \sum_{\mbf{X^J}\in\mc{H^I}} \s{\omega}^{\I\J} \left(\frac{|\sN{Y}^{IJ}|}{|\xi^{\I\J}|}\right)^{2m-2} \left(\frac{1}{|\xi^{\I\J}|}\right)^{2} \Delta V^{\J} . 
  \label{eqn:1d1}
\end{align}
Adopting a uniform discrete representation with $\Delta V^J = A \Delta X$ and $\delta = N \Delta X$ and the specific choice of influence function in \cref{eqn:omega}, the shape tensor is computed to be $L = 2 A N \Delta X$. Using a homogeneous deformation with $x = (1+a)X$,
\begin{equation}
  \sum_{\mbf{X^J}\in\mc{H^I}} \left(\s{\omega}^{\I\J}\right)^2 \left(\frac{|\sN{Y}^{IJ}|}{|\xi^{\I\J}|}\right)^{4m-2} \left(\frac{1}{|\xi^{\I\J}|}\right)^{2} \Delta V^{\J} = 2 \frac{A}{\Delta X} (1+a)^{4m-2} \sum_{p=1}^{N} \frac{1}{p^{2}} ,
  \label{eqn:1d2}
\end{equation}
\begin{equation}
  \sum_{\mbf{X^J}\in\mc{H^I}} \s{\omega}^{\I\J} \left(\frac{|Y^{IJ}|}{|\xi^{\I\J}|}\right)^{2m-2} \left(\frac{1}{|\xi^{\I\J}|}\right)^{2} \Delta V^{\J} = 2 \frac{A}{\Delta X} (1+a)^{2m-2} \sum_{p=1}^{N} \frac{1}{p^{2}} .
  \label{eqn:1d3}
\end{equation}
Using \cref{eqn:1d1,eqn:1d2,eqn:1d3},
\begin{equation}
  \p{f^\I}{x^\I} = - \frac{2 (1+a)^{2m-2}}{N (\Delta X)^2} \sum_{p=1}^{N} \frac{1}{p^{2}}\left(  \p{S_{(m)}}{\bar{E}_{(m)}} \frac{(1+a)^{2m}}{4N} + (2m-1) \, S_{(m)} \right) .
  \label{eqn:df1D}
\end{equation}

Next, the peridynamic Lagrangian Seth–Hill strain is computed to evaluate the stress state.
\begin{equation*}
  \bar{C}_{(m)} = (1+a)^{2m} , \quad \bar{E}_{(m)} = \frac{(1+a)^{2m} - 1}{2m} , \quad m \neq 0, 
\end{equation*}
and 
\begin{equation*}
  \bar{E}_{(m)} = \ln{(1+a)} , \quad m = 0.
\end{equation*}
A linear 1D constitutive model is used
\begin{equation}
  S_{(m)} = \kappa \bar{E}_{(m)} , \quad \p{S_{(m)}}{\bar{E}_{(m)}} = \kappa .
  \label{eqn:S1d}
\end{equation}
Substituting \cref{eqn:S1d} in \cref{eqn:df1D},
\begin{align*}
  \p{f^\I}{x^\I} = - 2\kappa \frac{(1+a)^{2m-2}}{N (\Delta X)^2} \sum_{p=1}^{N} \frac{1}{p^{2}} \bigg[ & \frac{(1+a)^{2m}}{4N} + \frac{2m-1}{2m} \Big( (1+a)^{2m} - 1 \Big) \bigg] , 
\end{align*}
for $m \neq 0$, and
\begin{equation*}
  \p{f^\I}{x^\I} = - 2\kappa \frac{(1+a)^{-2}}{N (\Delta X)^2} \sum_{p=1}^{N} \frac{1}{p^{2}} \bigg[ \frac{1}{4N} - \ln(1+a) \bigg] , \quad m = 0.
\end{equation*}
Stability enforces positivity of the terms inside the brackets. In other words, the model is unstable if the following conditions are violated: 
\begin{equation*}
  \bigg[ \frac{1}{4N} + \frac{2m-1}{2m} \bigg] (1+a)^{2m} > \frac{2m-1}{2m} , \quad m\neq0,
\end{equation*}
and 
\begin{equation*}
  \ln(1+a) < \frac{1}{4N} , \quad m=0.
\end{equation*}

\begin{figure}[!htbp]
  \centering
  \includegraphics[width=0.6\columnwidth]{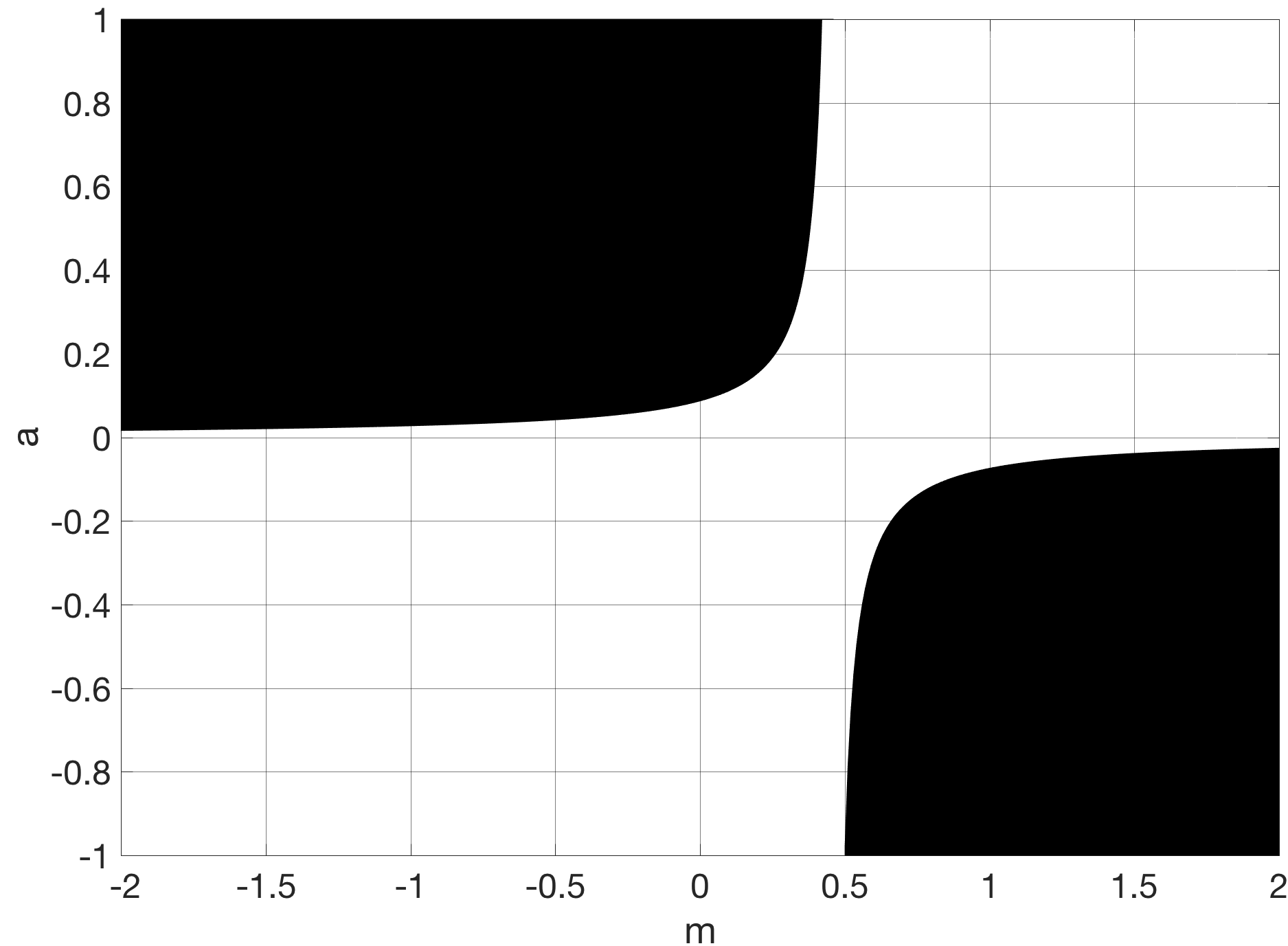}
  \caption{Stability of a discretized hyperelastic material with $N = 3$, under uniform hydrostatic deformation $\Dp{u_i}{X_j}=a\delta_{ij}$. Black color shows the unstable regions in 1D.}
  \label{fig:1d_hydro}
\end{figure}

Based on $m$ and $a$, stable and unstable regions are distinguished in \cref{fig:1d_hydro}, for a typical discretization with 3 nodes per horizon ($N=3$). As discussed in \cref{sec:hydrostatic}, compression is a problem for $m>1/2$, while instability is present in tensile mode for $m<1/2$.

\begin{remark}
\label{rmk:3}
  Increasing the nonlocality $N$ is unfavorable to the stability. Consider $m=1$, for instance, where $a_{\rm cr} = -0.184$ for $N=1$ and $a_{\rm cr} = -0.074$ for $N=3$ (the model is unstable for $a<a_{\rm cr}$). In other words, a greater compression is required to make the former unstable.
\end{remark}

\subsection{2D}
\label{sec:2D}
A discretized form of \cref{eqn:finaldf} in a two-dimensional setup is considered here.
\begin{align}
  \p{f^\I_k}{x^\I_l} = 
  - & \p{S_{(m)ij}}{\bar{E}_{(m)rs}} L^{-1}_{turs} L^{-1}_{pqij} \Delta V^\I \sum_{\mbf{X^J}\in\mc{H^I}} \Bigg[ \left(\s{\omega}^{\I\J}\right)^2 \bigg(\frac{\left|\s{Y}^{\mbf{IJ}}\right|}{\left|\bs{\xi^{\I\J}}\right|}\bigg)^{4m-2} \frac{\sN{Y}^{\I\J}_k \sN{Y}_l^{\I\J}}{\left|\s{Y}^{\mbf{IJ}}\right|^2} \frac{\xi^{\I\J}_p\xi^{\I\J}_q \xi_t^{\I\J} \xi_u^{\I\J}}{\left|\bs{\xi^{\I\J}}\right|^{6}} \Delta V^{\J} \Bigg] \nonumber \\
  & - 2 S_{(m)ij} L^{-1}_{pqij} \sum_{\mbf{X^J}\in\mc{H^I}} \Bigg\{ \s{\omega}^{\I\J} \bigg(\frac{\left|\s{Y}^{\mbf{IJ}}\right|}{\left|\bs{\xi^{\I\J}}\right|}\bigg)^{2m-2} \frac{\xi^{\I\J}_p\xi^{\I\J}_q}{\left|\bs{\xi^{\I\J}}\right|^{4}} \bigg[ 2(m-1)\frac{\sN{Y}_k^{\I\J} \sN{Y}^{\I\J}_l}{\left|\s{Y}^{\mbf{IJ}}\right|^2} + \delta_{kl} \bigg] \, \Delta V^{\J} \Bigg\} . 
  \label{eqn:2Dd1}
\end{align}

Considering a uniform discretization with $\Delta V^J = b (\Delta X)^2$ and $\delta = N \Delta X$ and an influence function as in \cref{eqn:omega}, the inverse of the shape tensor is computed first (see \cref{sec:AppendixA}). Then, \cref{eqn:2Dd1} is simplified for a general 2D uniform deformation with $x_i = X_i + a_{ij}X_j$.
\begin{align}
  \p{f^\I_k}{x^\I_l} = 
  - & \frac{1}{\pi^2N^2} \p{S_{(m)ij}}{\bar{E}_{(m)rs}} \sum_{\mbf{X^J}\in\mc{H^I}} \Bigg[ \bigg(\frac{\left|\s{Y}^{\mbf{IJ}}\right|}{\left|\bs{\xi^{\I\J}}\right|}\bigg)^{4m-2} \frac{\sN{Y}^{\I\J}_k \sN{Y}_l^{\I\J}}{\left|\s{Y}^{\mbf{IJ}}\right|^2} \frac{(4 \, \xi^{\I\J}_i\xi^{\I\J}_j - \xi^{\I\J}_p\xi^{\I\J}_p\delta_{ij}) (4 \, \xi^{\I\J}_r\xi^{\I\J}_s - \xi^{\I\J}_q\xi^{\I\J}_q\delta_{rs}) }{\left|\bs{\xi^{\I\J}}\right|^{6}} \Bigg] \notag \\
  & - \frac{2}{\pi N} S_{(m)ij} \sum_{\mbf{X^J}\in\mc{H^I}} \Bigg\{ \bigg(\frac{\left|\s{Y}^{\mbf{IJ}}\right|}{\left|\bs{\xi^{\I\J}}\right|}\bigg)^{2m-2} \frac{(4 \, \xi^{\I\J}_i\xi^{\I\J}_j - \xi^{\I\J}_n\xi^{\I\J}_n\delta_{ij})}{\left|\bs{\xi^{\I\J}}\right|^{4}} \bigg[ 2(m-1)\frac{\sN{Y}_k^{\I\J} \sN{Y}^{\I\J}_l}{\left|\s{Y}^{\mbf{IJ}}\right|^2} + \delta_{kl} \bigg] \Bigg\} ,
  \label{eqn:2Ddf}
\end{align}
where $Y_i^{IJ} = \xi_i^{IJ} + a_{ij}\xi_j^{IJ}$. It will be tedious, if possible, to analyze a general deformation with all the 4 $a_{ij}$ constants. 
Let consider a pure hydrostatic loading condition with $a_{11}=a_{22}=a$ and $a_{12}=a_{21}=0$. \cref{eqn:2Ddf} is then simplified to 
\begin{align}
  \p{f^\I_k}{x^\I_l} = 
  & - \frac{1}{\pi^2N^2} \p{S_{(m)ij}}{\bar{E}_{(m)rs}} (1+a)^{4m-2} \sum_{\mbf{X^J}\in\mc{H^I}} \Bigg[ \frac{\xi^{\I\J}_k\xi^{\I\J}_l(4 \, \xi^{\I\J}_i\xi^{\I\J}_j - \xi^{\I\J}_p\xi^{\I\J}_p\delta_{ij}) (4 \, \xi^{\I\J}_r\xi^{\I\J}_s - \xi^{\I\J}_q\xi^{\I\J}_q\delta_{rs}) }{\left|\bs{\xi^{\I\J}}\right|^{8}} \Bigg] \nonumber \\
  & - \frac{2}{\pi N} S_{(m)ij} \, (1+a)^{2m-2} \sum_{\mbf{X^J}\in\mc{H^I}} \frac{(4\xi^{\I\J}_i\xi^{\I\J}_j - \xi^{\I\J}_n\xi^{\I\J}_n\delta_{ij})}{\left|\bs{\xi^{\I\J}}\right|^{4}} \Bigg[ 2(m-1)\frac{\xi_k^{\I\J} \xi^{\I\J}_l}{\left|\bs{\xi^{\I\J}}\right|^2} + \delta_{kl} \Bigg] .
  \label{eqn:2Ddfs}
\end{align}
Using $\bxi^{\I\J} = p\Delta X \hat{x}_1 + q\Delta X \hat{x}_2$ and involving symmetry,
\begin{align*}
  \sum_{\mbf{X^J}\in\mc{H^I}} \frac{(\xi^{\I\J}_1)^4}{\left|\bs{\xi^{\I\J}}\right|^{6}} \Delta V^J 
  &= 4 \sum_{p+q=1}^{p^2+q^2 \leq N^2} \frac{p^4}{(p^2+q^2)^{3}} \frac{\Delta V^J}{(\Delta X)^2} , \nonumber \\
  &= 4\frac{\Delta V^J}{(\Delta X)^2} + 4 \sum_{p+q=2}^{p^2+q^2 \leq N^2} \frac{p^4}{(p^2+q^2)^{3}} \frac{\Delta V^J}{(\Delta X)^2} , \nonumber \\
  &\approx 4b + 4 \int_{\Delta X}^{\delta} \int_{0}^{\pi/2} \frac{(r cos\theta)^4}{r^{6}} b r {\rm d}\theta {\rm d}r ,
\end{align*}
where the last summation is approximated by an integral to obtain a closed-form solution. Then,
\begin{align*}
  \sum_{\mbf{X^J}\in\mc{H^I}} \frac{(\xi^{\I\J}_1)^4}{\left|\bs{\xi^{\I\J}}\right|^{6}} \Delta V^J  
  &= 4b + \frac{3\pi}{4} b \ \ln(\frac{\delta}{\Delta X}) , \nonumber \\
  &= b \left(4 + \frac{3\pi}{4} \ln N\right) .
\end{align*}
Similarly, the following summations are approximated using integration:
\begin{equation*}
  \sum_{\mbf{X^J}\in\mc{H^I}} \frac{1}{\left|\bs{\xi^{\I\J}}\right|^{2}} \Delta V^J = b \left(8 + 2\pi \ln N\right) ,
\end{equation*}
\begin{equation*}
  \sum_{\mbf{X^J}\in\mc{H^I}} \frac{(\xi^{\I\J}_1)^2}{\left|\bs{\xi^{\I\J}}\right|^{4}} \Delta V^J = b \left(4 + \pi \ln N\right) ,
\end{equation*}
\begin{equation*}
  \sum_{\mbf{X^J}\in\mc{H^I}} \frac{(\xi^{\I\J}_1)^2(\xi^{\I\J}_2)^2}{\left|\bs{\xi^{\I\J}}\right|^{6}} \Delta V^J = b \frac{\pi}{4} \ln N ,
\end{equation*}
\begin{equation*}
  \sum_{\mbf{X^J}\in\mc{H^I}} \frac{(\xi^{\I\J}_1)^6}{\left|\bs{\xi^{\I\J}}\right|^{8}} \Delta V^J = b \left(4 + \frac{5\pi}{8} \ln N\right) ,
\end{equation*}
\begin{equation*}
  \sum_{\mbf{X^J}\in\mc{H^I}} \frac{(\xi^{\I\J}_1)^4(\xi^{\I\J}_2)^2}{\left|\bs{\xi^{\I\J}}\right|^{8}} \Delta V^J = b \frac{\pi}{8} \ln N .
\end{equation*}
Next, by using symmetry and noting that only the even combinations of $\xi^{\I\J}_i$ produce nonzero values, the summations in \cref{eqn:2Ddfs} are simplified.
\begin{align*}
  \sum_{\mbf{X^J}\in\mc{H^I}} \frac{\xi^{\I\J}_i \xi^{\I\J}_j}{\left|\bs{\xi^{\I\J}}\right|^{4}} = \delta_{ij} \sum_{\mbf{X^J}\in\mc{H^I}} \frac{(\xi^{\I\J}_1)^2}{|\bs{\xi^{\I\J}}|^4} ,
\end{align*}
\begin{align*}
  & \sum_{\mbf{X^J}\in\mc{H^I}} \frac{\xi^{\I\J}_i \xi^{\I\J}_j \xi^{\I\J}_k \xi^{\I\J}_l}{\left|\bs{\xi^{\I\J}}\right|^{6}} = \left[ \delta_{ij}\delta_{kl}(1-\delta_{ik}) + \delta_{ik}\delta_{jl}(1-\delta_{ij}) + \delta_{il}\delta_{jk}(1-\delta_{ij}) \right] \sum_{\mbf{X^J}\in\mc{H^I}} \frac{(\xi^{\I\J}_1)^2 (\xi^{\I\J}_2)^2}{|\bs{\xi^{\I\J}}|^6} + \delta_{ij}\delta_{jk}\delta_{kl} \sum_{\mbf{X^J}\in\mc{H^I}} \frac{(\xi^{\I\J}_1)^4}{|\bs{\xi^{\I\J}}|^6} ,
\end{align*}
\begin{align*}
  \sum_{\mbf{X^J}\in\mc{H^I}} \frac{\xi^{\I\J}_i \xi^{\I\J}_j \xi^{\I\J}_k \xi^{\I\J}_l \xi^{\I\J}_r \xi^{\I\J}_s}{\left|\bs{\xi^{\I\J}}\right|^{8}} & = \gamma_{ijklrs} \sum_{\mbf{X^J}\in\mc{H^I}} \frac{(\xi^{\I\J}_1)^4 (\xi^{\I\J}_2)^2}{|\bs{\xi^{\I\J}}|^8} + \delta_{ij}\delta_{jk}\delta_{kl}\delta_{lr}\delta_{rs} \sum_{\mbf{X^J}\in\mc{H^I}} \frac{(\xi^{\I\J}_1)^6}{|\bs{\xi^{\I\J}}|^8} , 
\end{align*}
where $\gamma_{ijklrs}$ is the following:
\begin{align}
 \gamma_{ijklrs} 
  = & \delta_{ij}\delta_{jk}\delta_{kl}\delta_{rs}(1-\delta_{is}) + \delta_{ij}\delta_{jk}\delta_{kr}\delta_{ls}(1-\delta_{is}) + \delta_{ij}\delta_{jl}\delta_{lr}\delta_{ks}(1-\delta_{is}) + \delta_{ik}\delta_{kl}\delta_{lr}\delta_{js}(1-\delta_{is}) + \delta_{jk}\delta_{kl}\delta_{lr}\delta_{is}(1-\delta_{ij}) \nonumber \\
  & + \delta_{ij}\delta_{jk}\delta_{ks}\delta_{lr}(1-\delta_{ir}) + \delta_{ij}\delta_{jl}\delta_{ls}\delta_{kr}(1-\delta_{ir}) + \delta_{ik}\delta_{kl}\delta_{ls}\delta_{jr}(1-\delta_{ir}) + \delta_{jk}\delta_{kl}\delta_{ls}\delta_{ir}(1-\delta_{ij}) + \delta_{ij}\delta_{jr}\delta_{rs}\delta_{kl}(1-\delta_{il}) \nonumber \\
  & + \delta_{ik}\delta_{kr}\delta_{rs}\delta_{jl}(1-\delta_{il}) + \delta_{jk}\delta_{kr}\delta_{rs}\delta_{il}(1-\delta_{ij}) + \delta_{il}\delta_{lr}\delta_{rs}\delta_{jk}(1-\delta_{ik}) + \delta_{jl}\delta_{lr}\delta_{rs}\delta_{ik}(1-\delta_{ij}) + \delta_{kl}\delta_{lr}\delta_{rs}\delta_{ij}(1-\delta_{is}) .
  \label{eqn:gamma}
\end{align}

Using the above calculations to simplify \cref{eqn:2Ddfs} (detailed in \cref{sec:appendixB}),
\begin{align}
  \p{f^\I_k}{x^\I_l} = & 
  - \frac{4}{\pi^2 N^2 (\Delta X)^2} (1+a)^{4m-2} \p{S_{(m)ij}}{\bar{E}_{(m)rs}} \Bigg\{ \gamma_{ijklrs} \frac{\pi}{2} \ln N + \delta_{ki}\delta_{ij}\delta_{jr}\delta_{rs}\delta_{sl} \Big( 16 + \frac{5\pi}{2} \ln N \Big) + \delta_{ij}\delta_{rs}\delta_{kl} \Big( 1 + \frac{\pi}{4} \ln N \Big) \nonumber \\
  & \qquad \qquad \qquad \qquad \qquad \qquad \quad \ \ - \delta_{rs} \bigg[ \Big( \delta_{ij}\delta_{kl}(1-\delta_{ik}) + \delta_{ik}\delta_{jl}(1-\delta_{ij}) + \delta_{il}\delta_{jk}(1-\delta_{ij}) \Big) \frac{\pi}{4}\ln N + \delta_{ki}\delta_{ij}\delta_{jl} \Big(4 + \frac{3\pi}{4}\ln N\Big) \bigg] \nonumber \\
  & \qquad \qquad \qquad \qquad \qquad \qquad \quad \ \ - \delta_{ij} \bigg[ \Big( \delta_{rs}\delta_{kl}(1-\delta_{rk}) + \delta_{rk}\delta_{sl}(1-\delta_{rs}) + \delta_{rl}\delta_{sk}(1-\delta_{rs}) \Big) \frac{\pi}{4}\ln N + \delta_{kr}\delta_{rs}\delta_{sl} \Big(4 + \frac{3\pi}{4}\ln N\Big) \bigg] \Bigg\} \nonumber \\
  & - \frac{4}{\pi N (\Delta X)^2} (1+a)^{2m-2} S_{(m)ij} \Bigg\{ 4(m-1) \bigg[ \Big( \delta_{ij}\delta_{kl}(1-\delta_{ik}) + \delta_{ik}\delta_{jl}(1-\delta_{ij}) + \delta_{il}\delta_{jk}(1-\delta_{ij}) \Big) \frac{\pi}{4}\ln N + \delta_{ki}\delta_{ij}\delta_{jl} \Big(4 + \frac{3\pi}{4}\ln N\Big) \bigg] \nonumber \\
  & \qquad \qquad \qquad \qquad \qquad \qquad  + (2-m) \delta_{ij} \delta_{kl} \Big(4 + \pi \ln N\Big) \Bigg\} .
  \label{eqn:2Ddfss}
\end{align}

Because the deformations are limited to hydrostatic, a simple constitutive law is adopted as
\begin{equation}
  S_{(m)ij} = \kappa \bar{E}_{(m)ij}, \quad \p{S_{(m)ij}}{\bar{E}_{(m)rs}} = \kappa \delta_{ij}\delta_{rs} ,
  \label{eqn:hydro2Ds}
\end{equation}
where $\kappa$ is twice the bulk modulus in a two-dimensional problem. Note that there is no shear for hydrostatic loadings. 

For a pure homogeneous hydrostatic loading with $x = (1+a)X$
\begin{equation}
  \bar{C}_{(m)ij} = (1+a)^{2m}\delta_{ij}, \quad \bar{E}_{(m)ij} = \frac{(1+a)^{2m} - 1}{2m}\delta_{ij} , \quad m\neq0, 
  \label{eqn:hydro2De}
\end{equation}
and
\begin{equation*}
  \bar{E}_{(m)ij} = \ln(1+a)\delta_{ij} , \quad m=0.
\end{equation*}
Substitute these back in \cref{eqn:2Ddfss} and simplify (see \cref{sec:appendixC}):
\begin{align*}
  \p{f^\I_k}{y^\I_l} 
  &= - \frac{4 (4 + \pi \ln N)}{\pi^2 N^2 (\Delta X)^2} \kappa (1+a)^{4m-2} \left[ 1 + 2 \pi N \Big( 1 - (1+a)^{-2m} \Big) \right] \delta_{kl} ,
\end{align*}
for $m \neq 0$, and
\begin{align*}
  \p{f^\I_k}{y^\I_l} 
  &= - \frac{4 (4 + \pi \ln N)}{\pi^2 N^2 (\Delta X)^2} \kappa (1+a)^{-2} \left[ 1 + 2 \pi N \ln(1+a) \right] \delta_{kl} , \quad m=0. 
\end{align*}
The eigenvalues of this diagonal matrix are identified as:
\begin{align*}
  \lambda_i &= - \frac{4 (4 + \pi \ln N)}{\pi^2 N^2 (\Delta X)^2} \kappa (1+a)^{4m-2} \left[ 1 + 2 \pi N \Big( 1 - (1+a)^{-2m} \Big) \right] , \quad i=1,2 
\end{align*}
for $m \neq 0$, and
\begin{align*}
  \lambda_i &= - \frac{4 (4 + \pi \ln N)}{\pi^2 N^2 (\Delta X)^2} \kappa (1+a)^{-2} \left[ 1 + 2 \pi N \ln(1+a) \right] , \quad i=1,2 , \quad m=0.
\end{align*}

\begin{figure}[!t]
  \centering
  \includegraphics[width=0.5\columnwidth]{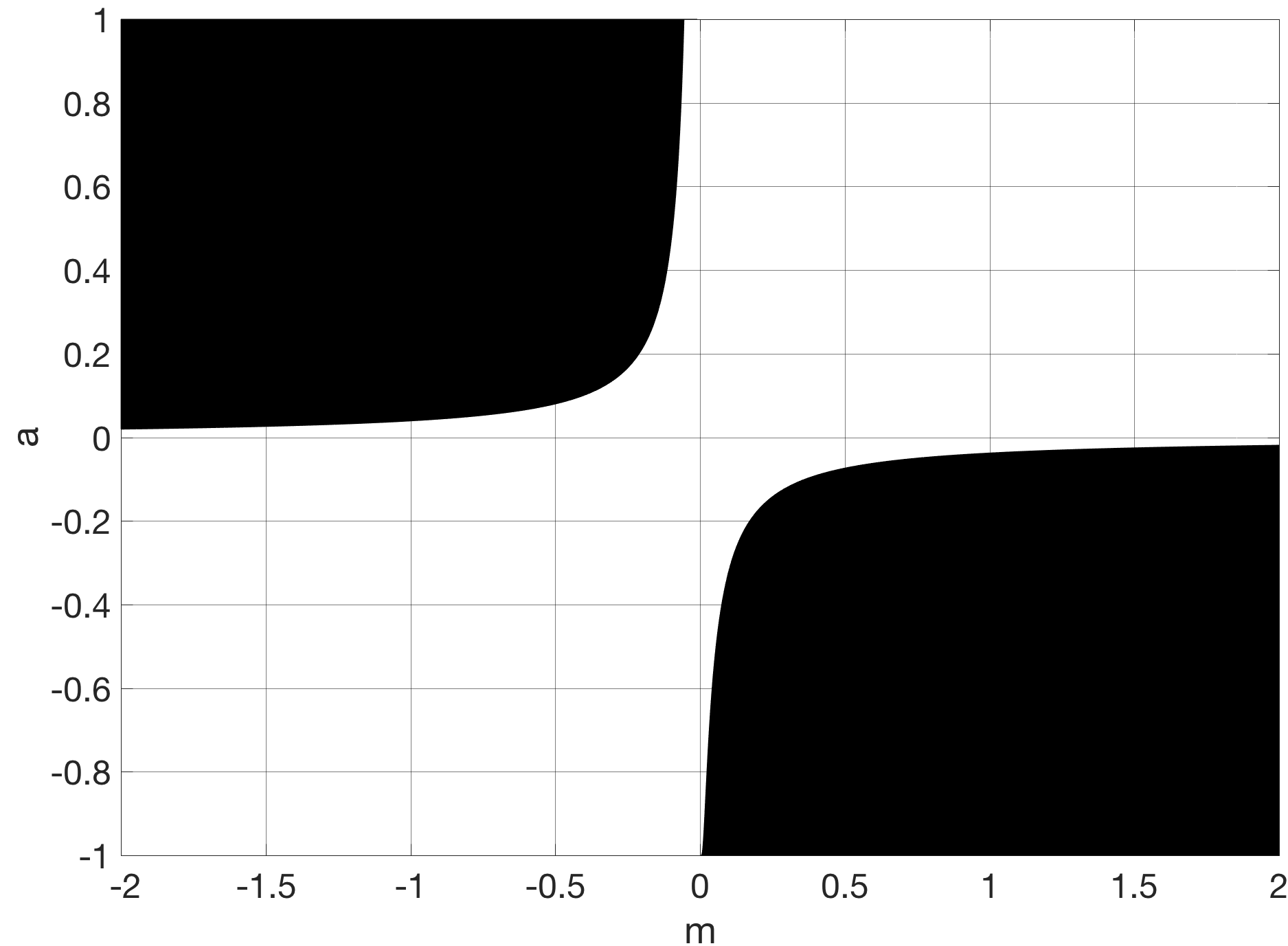}
  \caption{Stability of a discretized hyperelastic material with $N = 3$, under uniform hydrostatic deformation $\Dp{u_i}{X_j}=a\delta_{ij}$. Black color shows the unstable regions in 2D.}
  \label{fig:2d_hydro}
\end{figure}

An equilibrium point is unstable if it does not satisfy the following criteria:
\begin{align*}
  1 + 2 \pi N \Big[ 1 - (1+a)^{-2m} \Big] < 0 , \quad m \neq 0 ,
\end{align*}
\begin{align*}
  1 + 2 \pi N \ln(1+a) < 0 , \quad m = 0 .
\end{align*}
Similar to the 1D case, the stable and unstable loading conditions are recognized for different $m$'s, depicted in \cref{fig:2d_hydro}

\begin{remark}
\label{rmk:4}
Increasing the dimension expands the unstable regions. For example, consider $m=1$ and $N=3$ in 2D, where $a_{cr}=-0.037$. This critical value is half the 1D case (cf. \cref{rmk:3}), meaning that a lower compression would cause instability. It can be shown that the region of instability is also larger in 3D than 2D. 
\end{remark}

\begin{remark}
\label{rmk:5}
Comparing the results of the discrete level from \cref{sec:discrete} with the continuum scale from \cref{sec:stability}, the stable regions are evidently larger in the discrete form. In other words, numerical representation helps the Lyapunov stability in the case of homogeneous hydrostatic deformations. 
\end{remark}

\begin{remark}
\label{rmk:6}
For $m=0$ and a small window of $m<0$ (in the discrete form), the model may be stable for any uniform hydrostatic deformation. However, other types of loading conditions may impose instability, which actually is the case and will be discussed in \cref{sec:m0stability}.
\end{remark}

While obtaining a closed-form solution for more complex loading scenarios is not easily achievable, stability of the model under such conditions can be examined through numerics. 

\section{Numerical examples}
\label{sec:numerical}
In this section, stability of the finite deformation model is numerically examined for different problems. The two special members of the material class, $m=0$ (Hencky strain) and $m=1$ (Green-Lagrange strain) are the primary targets.

\subsection{1D singular bar}
\label{sec:singular}
The problem was designed by \citet{breitenfeld2014non} to show the unstable behavior of the original correspondence theory. 
%
It involves a 1D bar characterized with a spatially varying Young's modulus
\begin{align*}
  E(x) =
  \begin{cases}
    \qquad \qquad \quad E_0 \qquad & x \leq L/2, \\
    E_0 \Big(1+\dfrac{\alpha}{2\sqrt{x/L-1/2}}\Big)^{-1} \quad & x > L/2, 
  \end{cases}
\end{align*}
fixed at one end and loaded with a stress $\sigma$ at the other end. The local analytical solution gives the normalized axial displacement 
\begin{align*}
  \dfrac{E_0 u(x)}{\sigma L} =
  \begin{cases}
    \qquad \qquad x/L \qquad & x \leq L/2, \\
    x/L + \alpha \sqrt{x/L-1/2} \quad & x > L/2, 
  \end{cases}
\end{align*}
and the normalized strain
\begin{align*}
  \dfrac{E_0 \eps(x)}{\sigma} =
  \begin{cases}
    \qquad \qquad 1 \qquad & x \leq L/2, \\
    1 + \dfrac{\alpha}{2\sqrt{x/L-1/2}} \quad & x > L/2, 
  \end{cases}
\end{align*}
 along the bar. 

\begin{figure*}[b!]
  \centering
  \subfloat[][Tension - displacement]{\includegraphics[width=0.47\textwidth,trim={0 1.93cm 0 0 },clip]{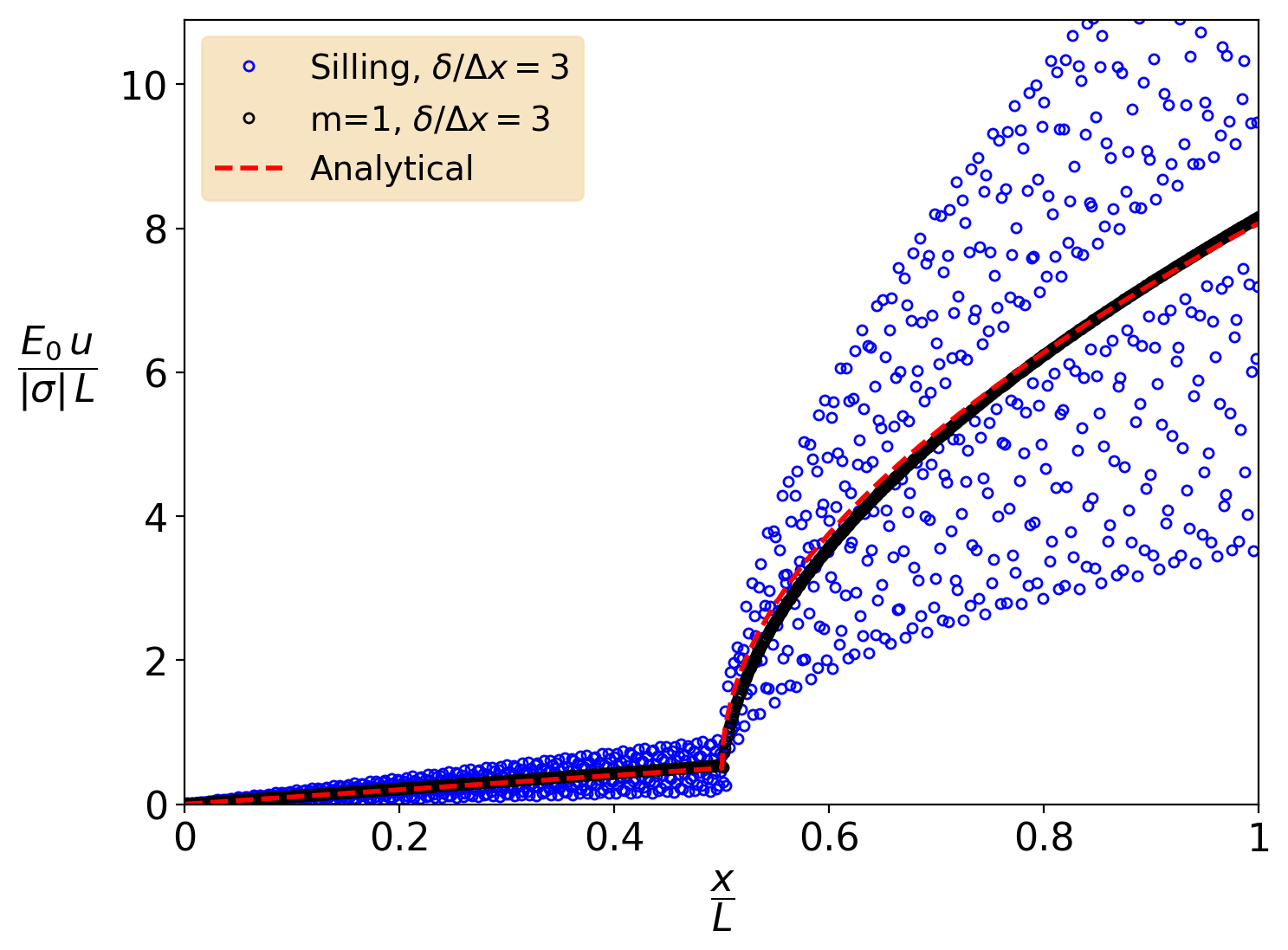}}
  \hspace*{0.7cm}
  \subfloat[][Tension - strain]{\includegraphics[width=0.47\textwidth,trim={0 1.93cm 0 0 },clip]{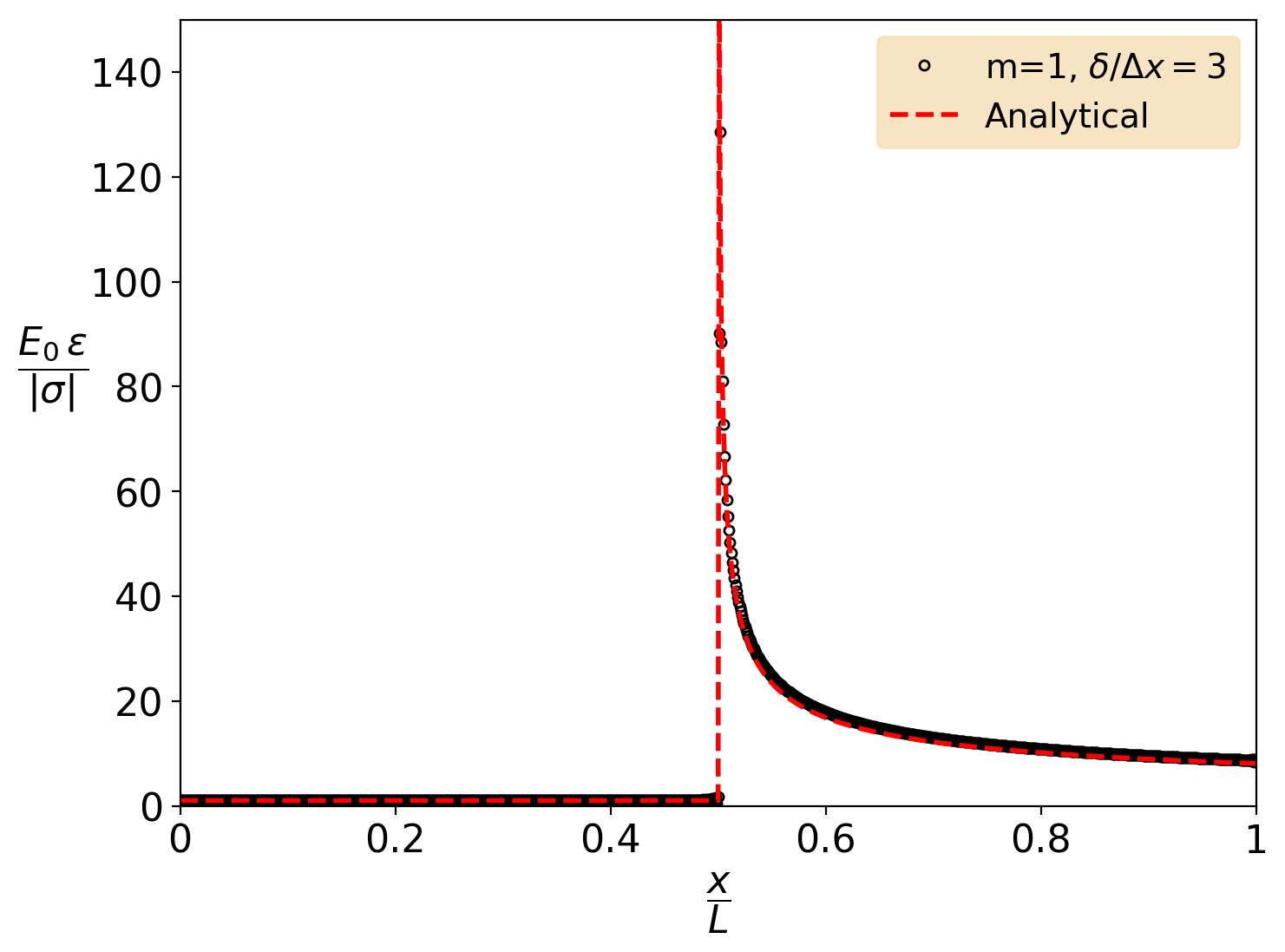}}

  \subfloat[][Compression - displacement]{\includegraphics[width=0.47\textwidth]{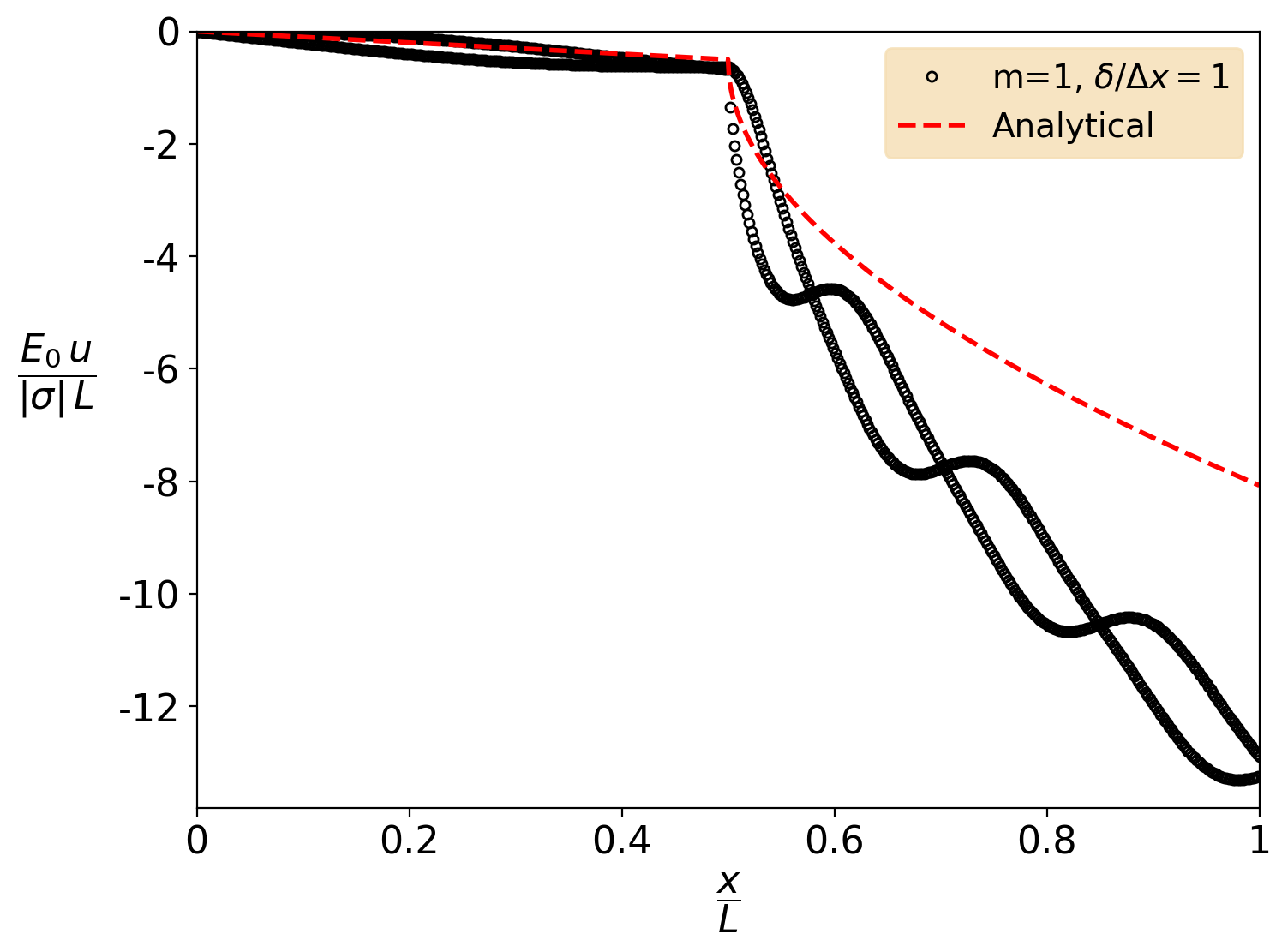}}
  \hspace*{0.7cm}
  \subfloat[][Compression - strain]{\includegraphics[width=0.47\textwidth]{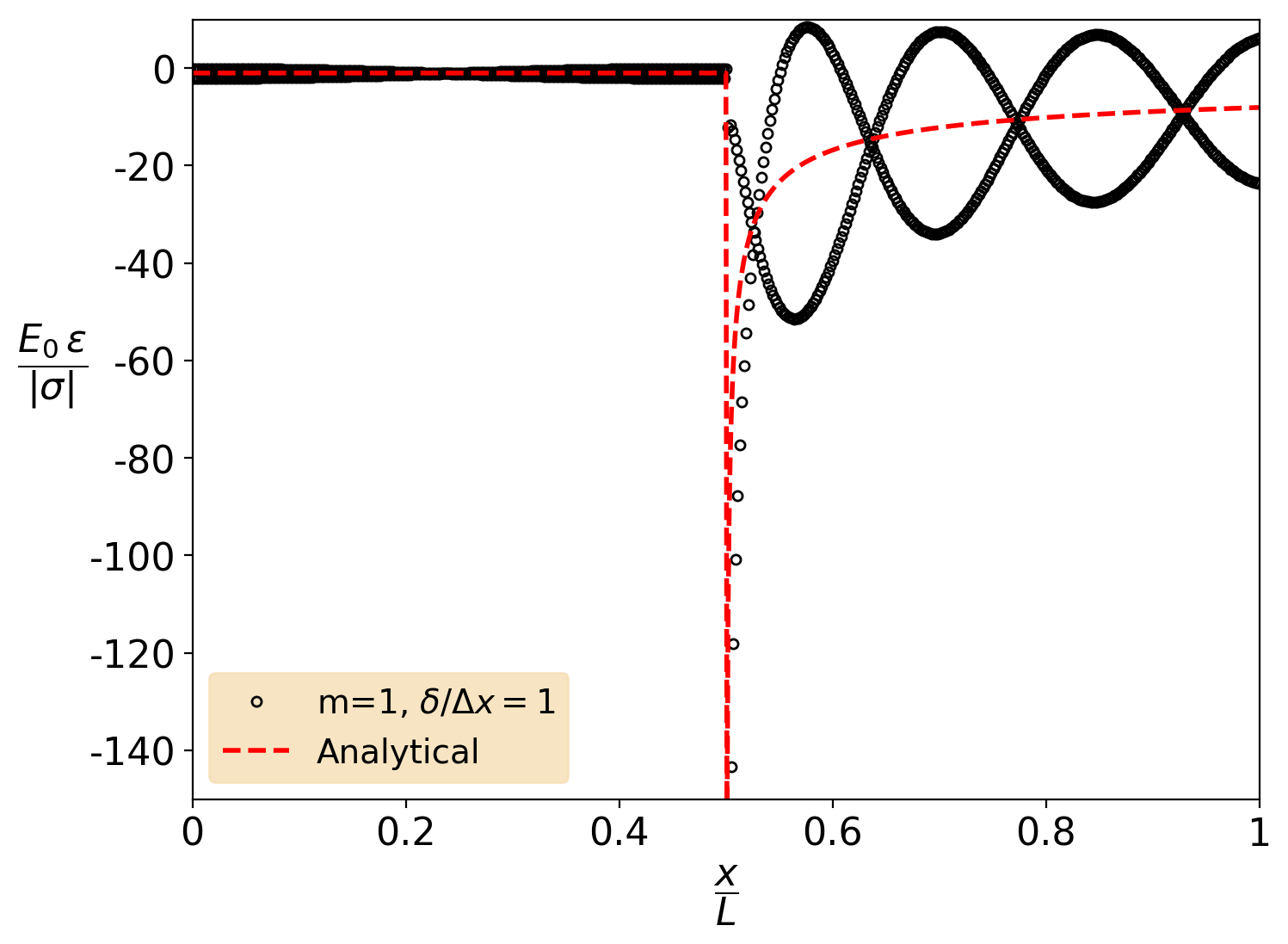}}
  \caption{While the $m=1$ model correctly predicts the response of the singular bar under tensile loading (a--b), it exhibits instability modes of deformation in dealing with compression (c--d).}
  \label{fig:singular_bar}
\end{figure*}

This static problem is solved using the $m=1$ model for an applied tensile stress (with $\sigma/E_0=10^{-3}$) and a smaller compressive stress (with $\sigma/E_0=-10^{-5}$). Normalized displacement and strain values are shown in \cref{fig:singular_bar}, when $\alpha = 10$. While the model correctly predicts the response of the bar with no apparent instability, the unstable behavior in compression mode leads to wrong results. The original correspondence model, in contrast, has a similar unstable behavior in both tensile and compressive modes \citep{breitenfeld2014non}.

\subsection{Influence of step size in different loading regimes}
\label{sec:stepsize}
\begin{figure*}[!ht]
  \centering
  \subfloat[][applied $\eps = 100 \% $]{\includegraphics[width=0.33\textwidth]{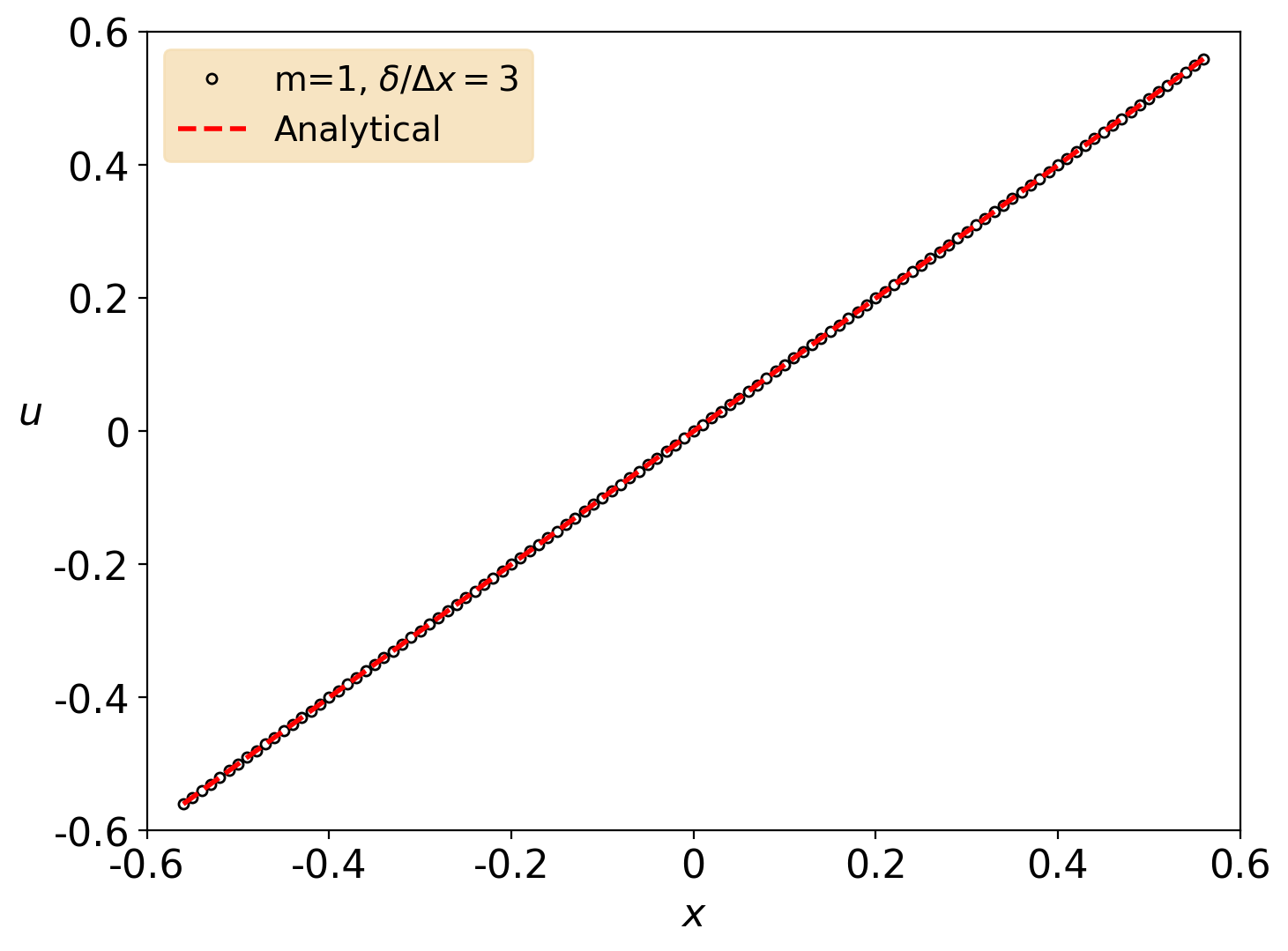}}
  \subfloat[][applied $\eps = -0.01 \% $]{\includegraphics[width=0.34\textwidth,trim={1cm 0 0 0},clip]{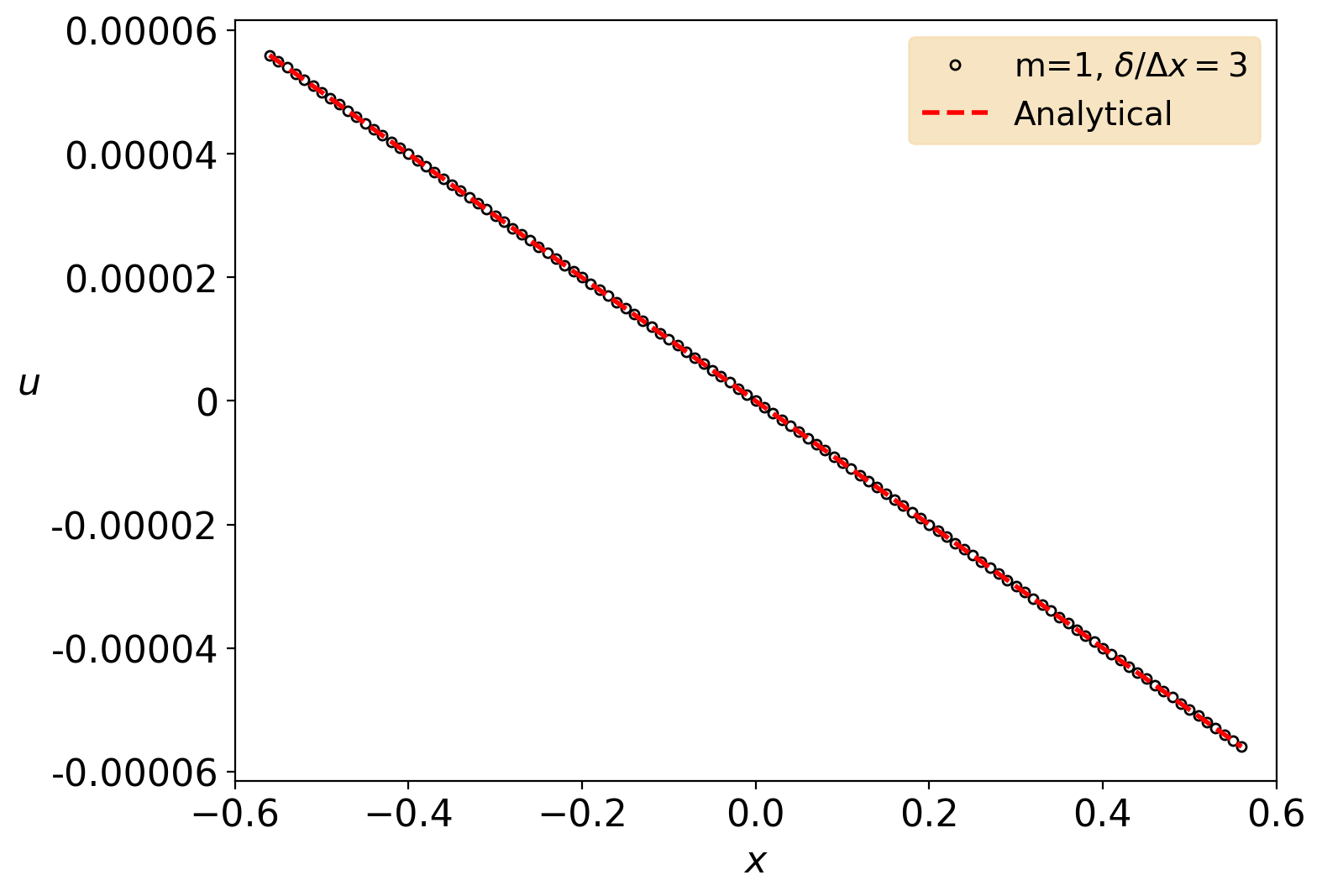}}
  \subfloat[][applied $\eps = -0.1 \% $]{\includegraphics[width=0.34\textwidth,trim={1cm 0 0 0},clip]{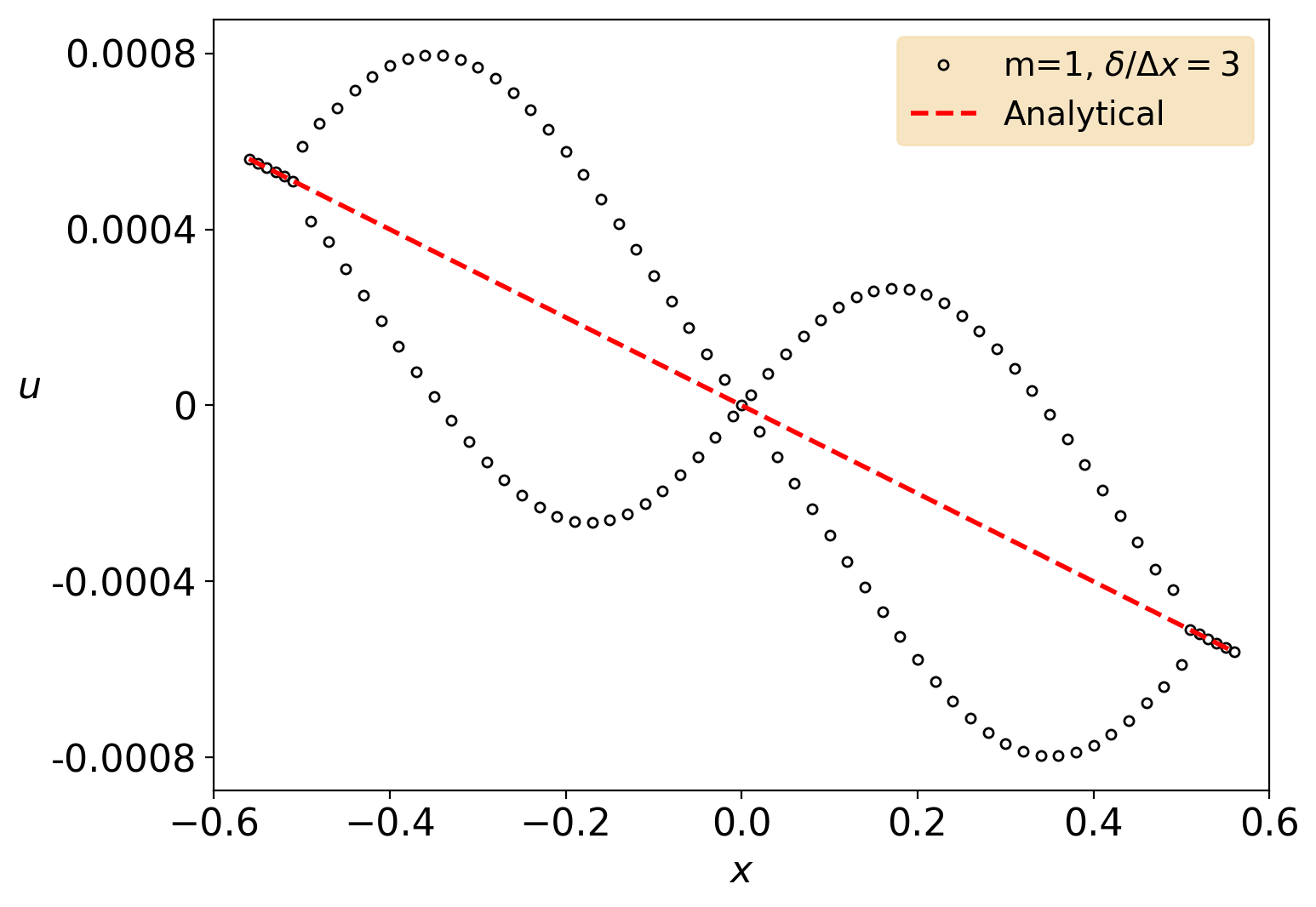}}
  \caption{Stability of the $m=1$ model in 1D setting. Significant differences in stable step size between compression and tension.}
  \label{fig:1d_stepSize}
\end{figure*}

The stable load step size can be significantly different between loading modes, even in the stable regions. For example, consider a 1-dimensional hyperelastic bar subjected to prescribed displacements at its ends. The equilibrium state is calculated using the $m=1$ material model with 3 nodes per horizon ($N=3$). Note that $a_{cr} = -7.4\%$ for this problem (cf. \cref{rmk:3}). The bar is initially undeformed and boundary displacements are applied instantly (in one step) over a boundary region that contains 6 nodes at each end. As shown in \cref{fig:1d_stepSize}, under tensile mode, the theory correctly finds the solution for an extremely large step size. However, the stable load step is several orders of magnitude smaller when compressive loading is applied, even when the compression state is well inside the stable region.

\subsection{Stability of the logarithmic model ($m=0$)}
\label{sec:m0stability}
\begin{figure*}[!h]
  \centering
  \subfloat[][Uniaxial tension]{\includegraphics[width=0.45\textwidth]{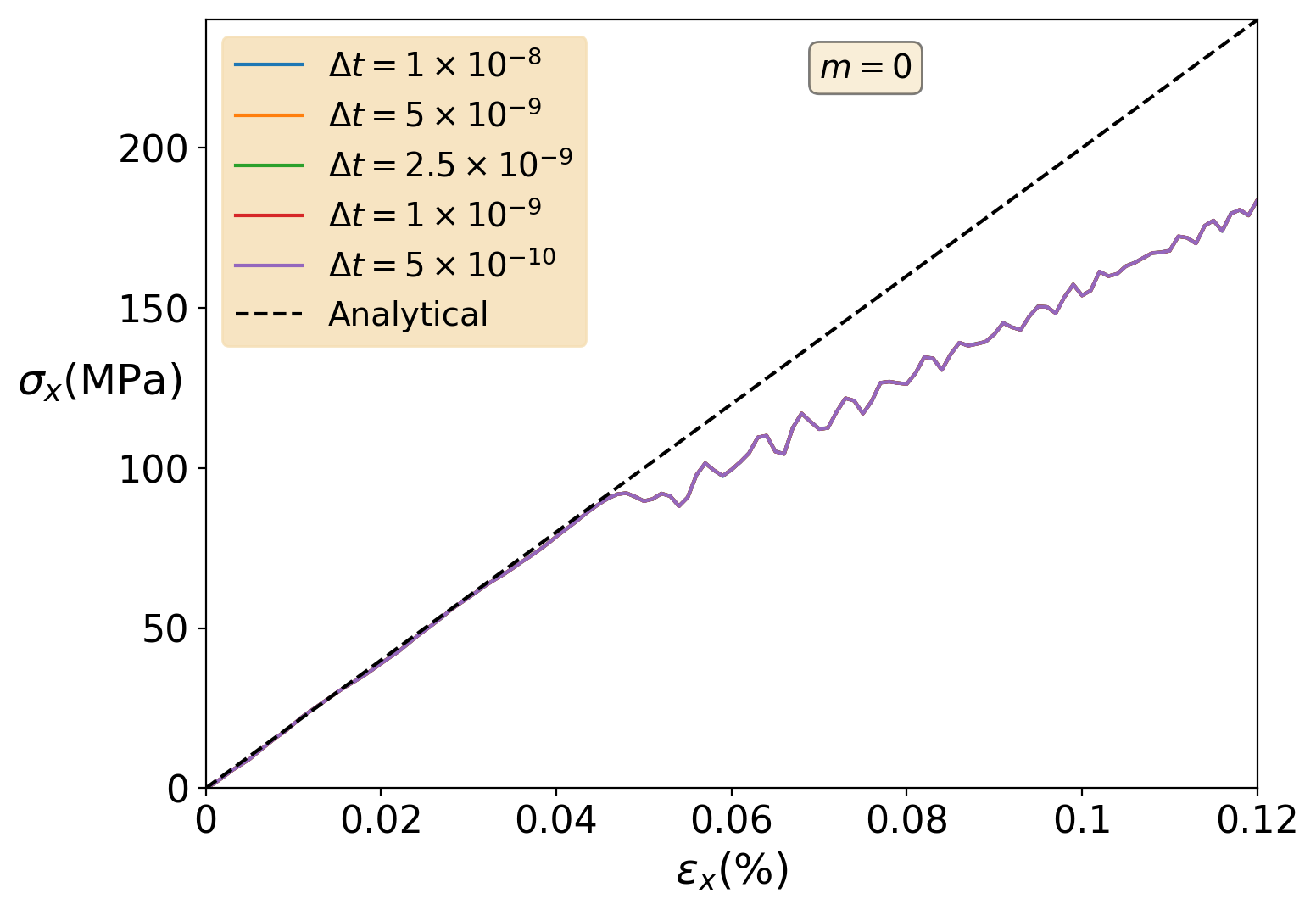}}
  \hspace*{0.5cm}
  \subfloat[][Uniaxial compression]{\includegraphics[width=0.45\textwidth]{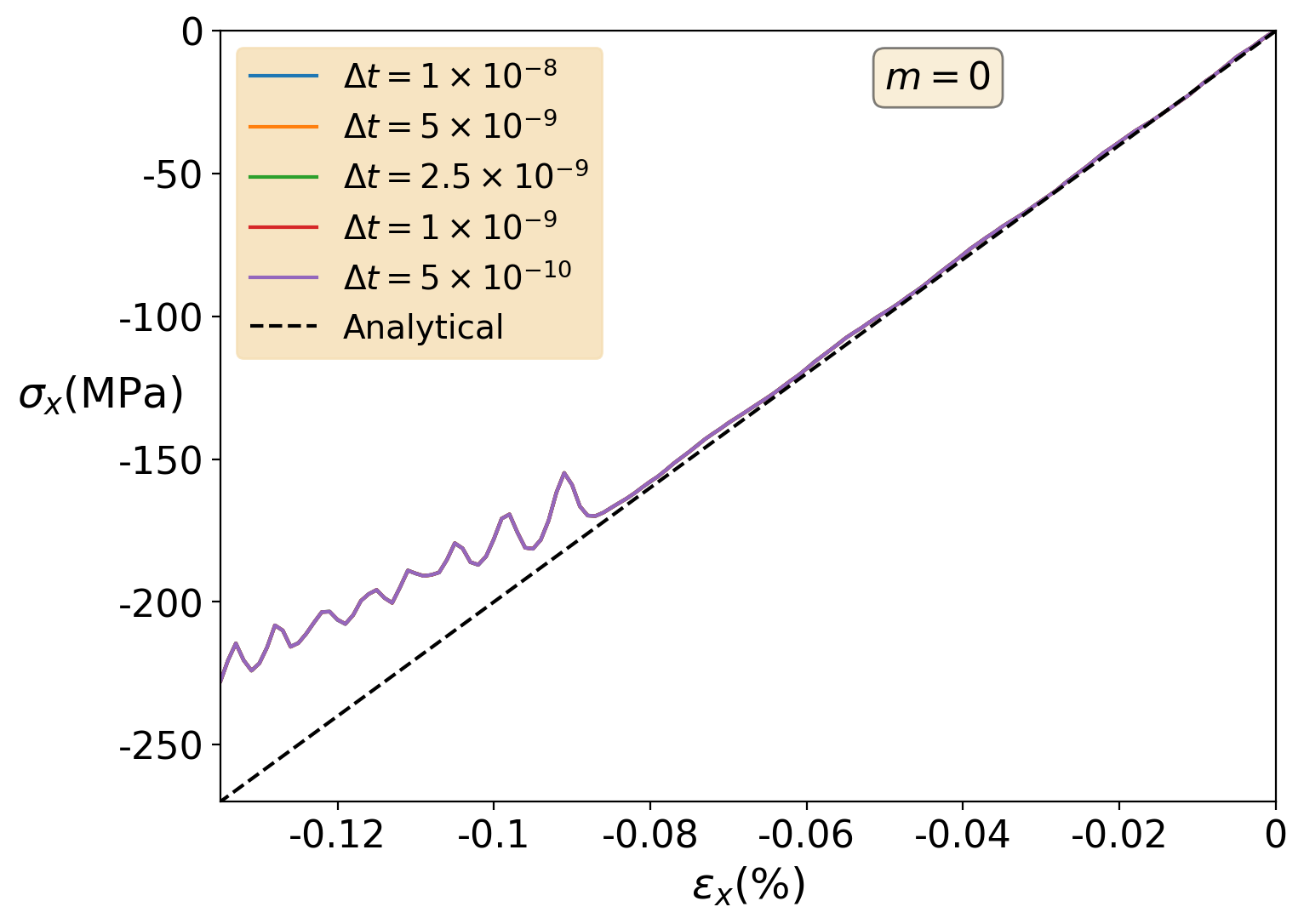}}

  \subfloat[][X displacement]{\includegraphics[width=0.32\textwidth]{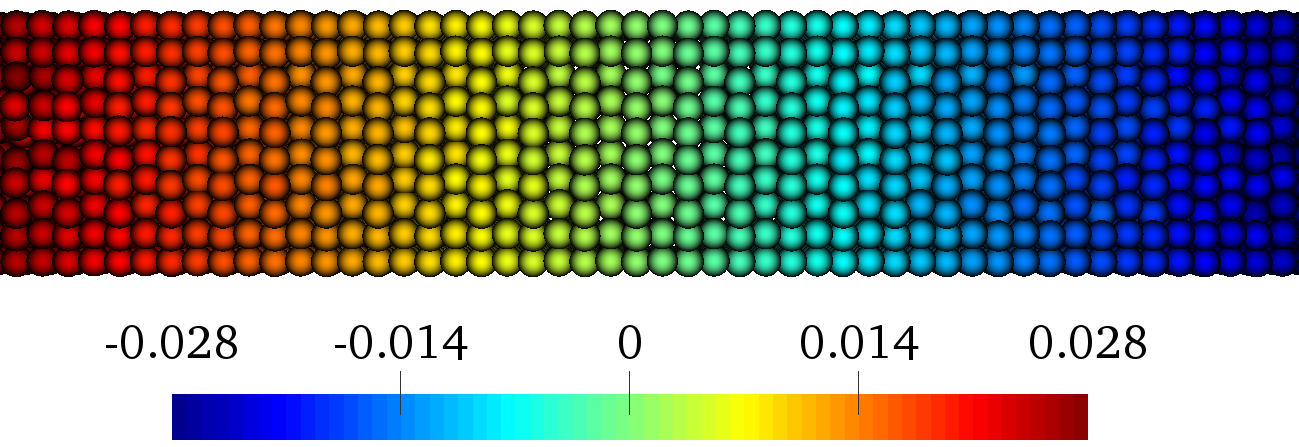}}
  \hspace*{0.25cm}
  \subfloat[][Y displacement]{\includegraphics[width=0.32\textwidth]{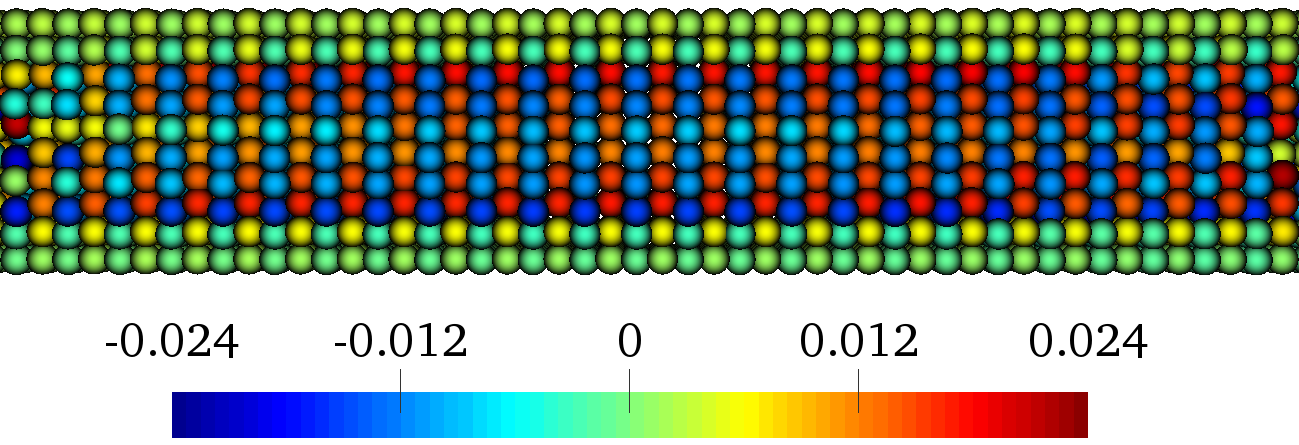}}
  \hspace*{0.25cm}
  \subfloat[][Z displacement]{\includegraphics[width=0.32\textwidth]{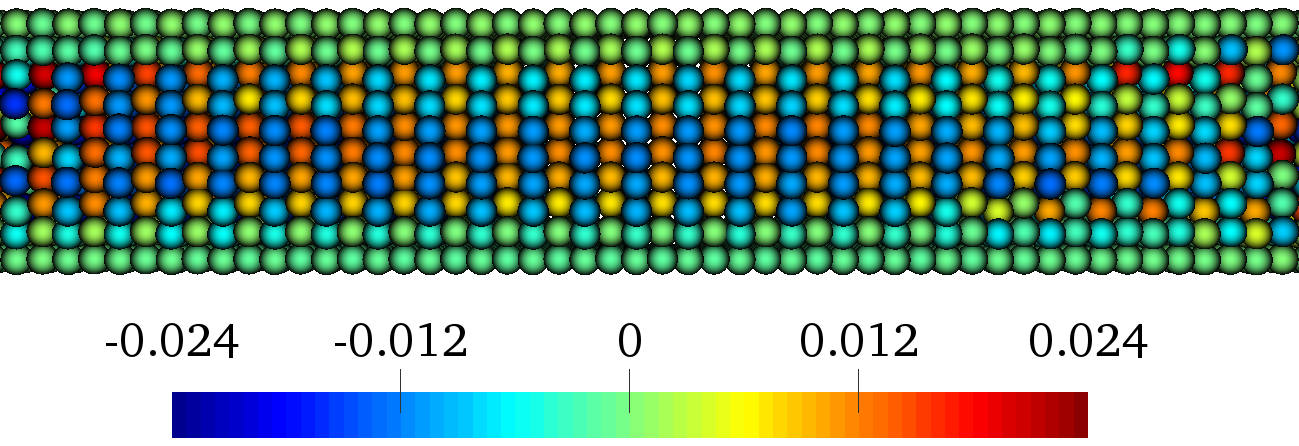}}
  \caption{Simulation results of a 3D hyperperelastic cuboid under uniaxial loading (along $X$-axis), using $m=0$ model. The predicted macroscopic response shows an unstable behavior under relatively small tensile (a) and compressive (b) loadings. While $u_x$ (c) appears as a smooth gradient field, instabilities are apparent in $u_y$ (d) and $u_z$ (e) values (displacements correspond to the compression case). Note that reducing the time step is not effective.}
  \label{fig:3d_uniaxial}
\end{figure*}

As discussed in \cref{sec:2D} and \cref{rmk:6}, it appears that the $m=0$ material model is stable in dealing with uniform hydrostatic loadings in 2D and 3D. However, other types of deformation can cause instability. Consider a simple problem involving a 3D hyperelastic cuboid ($4\times1\times1$) under uniaxial loading (along $X$-axis). Explicit time integration is employed to simulate the problem with prescribed displacements at the ends of the material. As shown in \cref{fig:3d_uniaxial}, the model becomes unstable under relatively small tension or compression states. It is seen that using a smaller time step has no effect on the predicted macroscopic behavior in this case. While the $X$-displacement values appear in a smooth gradient, instabilities are evident in the $Y$- and $Z$-displacements. 

\subsection{Waves in a hyperelastic bar}
\label{sec:waves}
As suggested by \citet{belytschko2000unified}, stability of a 1D hyperelastic bar is analyzed here. The material is modeled by different members of the finite deformation family, upon introduction of a perturbation deformation in the form of a planar wave 
\begin{equation}
  u = u_0 {\rm e}^{{\rm i}(kX-\omega t)} ,
  \label{eqn:dispersion0}
\end{equation}
in which $k$ is the wave number and $\omega$ is a complex angular frequency. 

\begin{figure*}[!htbp]
  \centering
  \subfloat[][]{\includegraphics[height=0.39\textwidth]{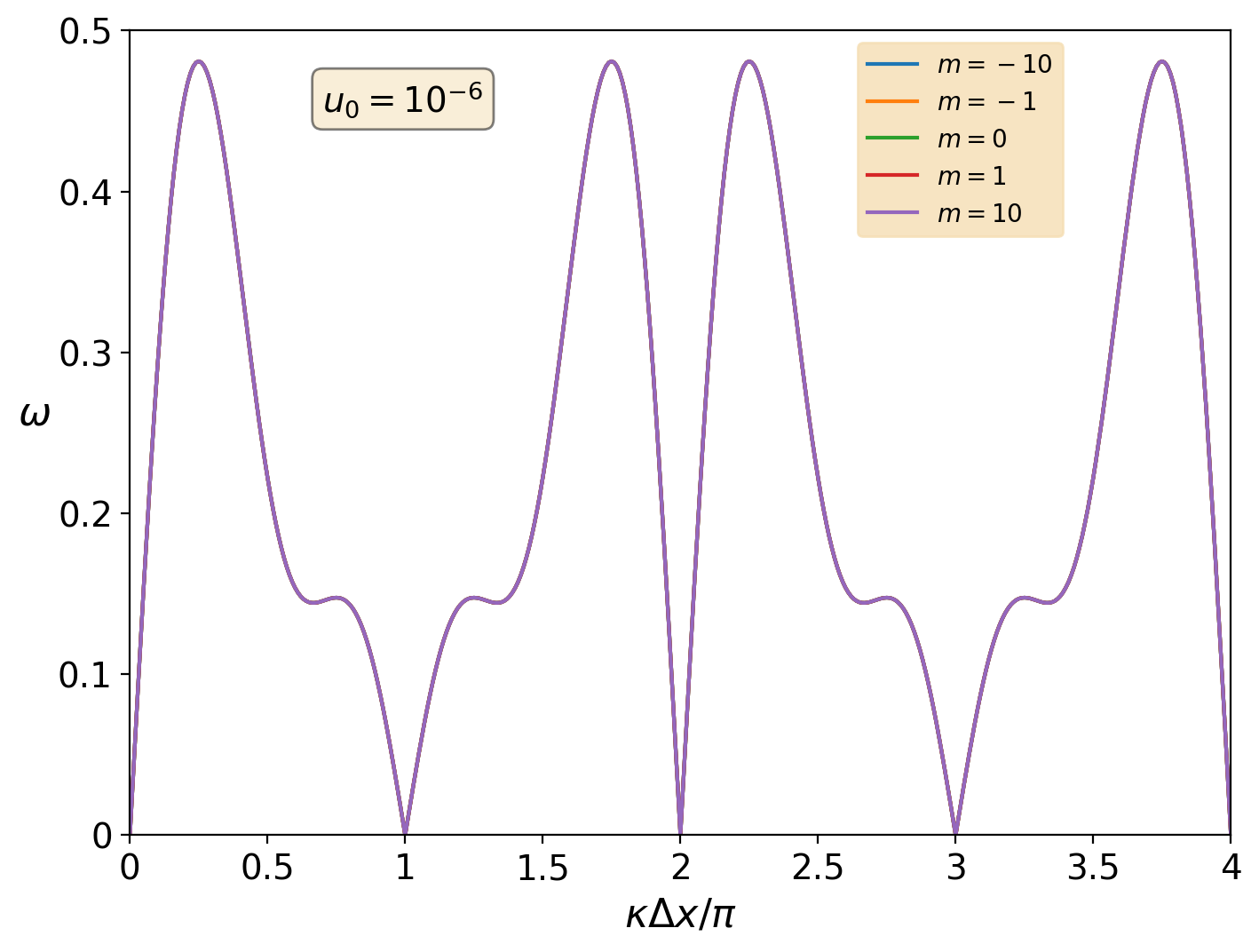}}
  \hspace*{0.4cm}
  \subfloat[][]{\includegraphics[height=0.39\textwidth,trim={1.7cm 0 0 0 },clip]{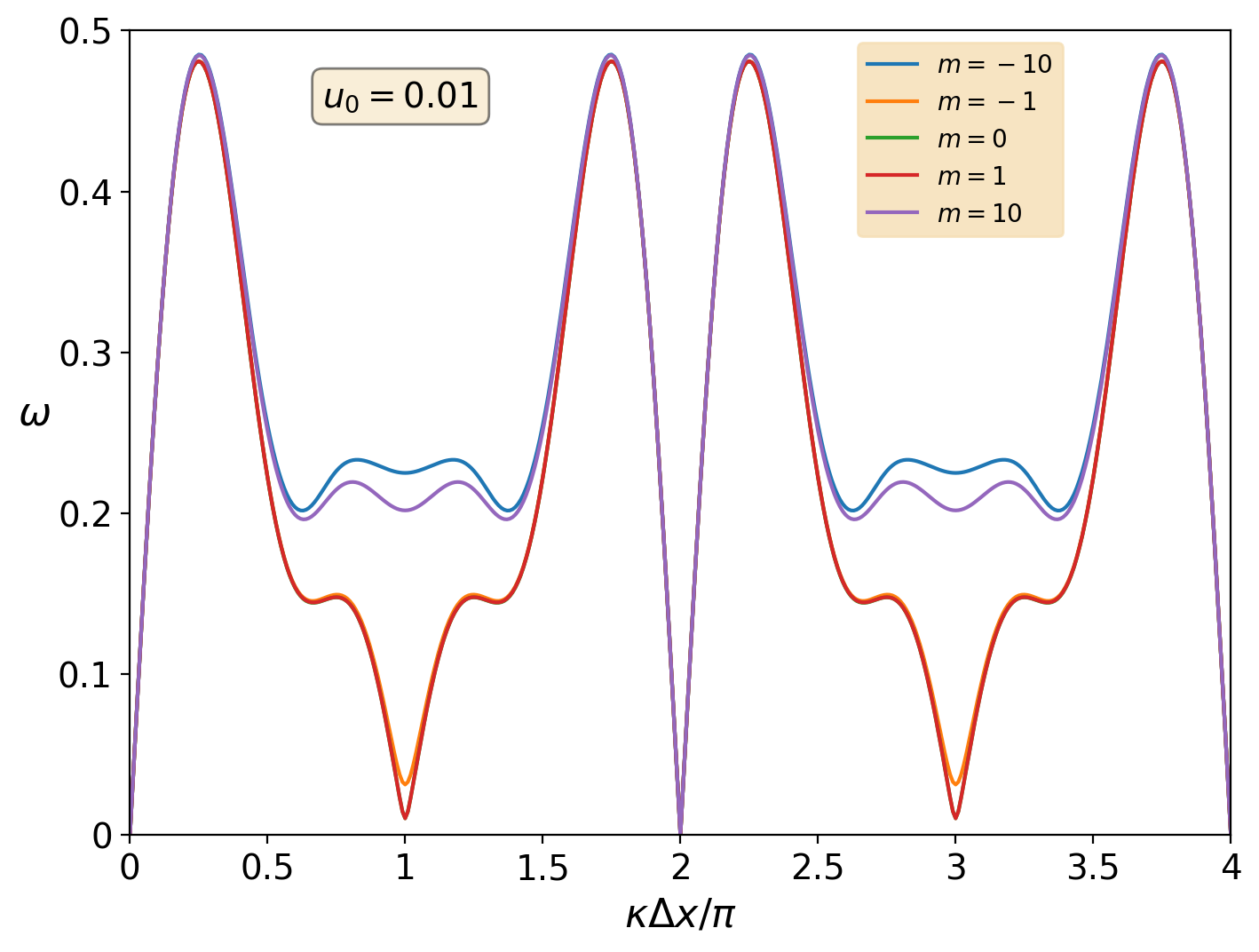}}
  \caption{Dispersion curves in a hyperelastic bar. Zero angular frequencies are indicative of instability. Behavior of different members of the peridynamic class are similar for small $u_0$ (a) and different for larger perturbations (b). All the models involve periodic instability points.}
  \label{fig:dispersion}
\end{figure*}

In 1D, the shape tensor $L = \int_{\mc{H}} {\rm d}\xi = \omega_0$. The force state in \cref{eqn:TC} then reads
\begin{equation}
  \sN{T}[\xi][X] = \frac{\s{\omega}[\xi]}{\omega_0} S_{(m)}(X) \, \bigg(\frac{\sN{Y}^2}{\xi^2}\bigg)^{m-1} \frac{\sN{Y}}{\xi^2}, 
  \label{eqn:dispersion1}
\end{equation}
and
\begin{equation}
  \sN{T}[-\xi][X+\xi] = \frac{\s{\omega}[\xi]}{\omega_0} S_{(m)}(X+\xi) \, \bigg(\frac{\sN{Y}^2}{\xi^2}\bigg)^{m-1} \frac{-\sN{Y}}{\xi^2}, 
  \label{eqn:dispersion2}
\end{equation}
where 
\begin{align}
  &\sN{Y}[\xi][X] = \xi + u_0 {\rm e}^{{\rm i}(kX-\omega t)} \Big({\rm e}^{{\rm i}k\xi}-1\Big) , \notag \\ 
  &= \xi + u_0 \Big[ \cos(kX-\omega t) + {\rm i} \sin(kX-\omega t) \Big] \Big[\cos(k\xi) + {\rm i} \sin(k\xi) - 1 \Big] .
  \label{eqn:dispersion3}
\end{align}
Using \cref{eqn:dispersion0,eqn:dispersion1,eqn:dispersion2}, the 1D peridynamic equation of motion is reduced as
\begin{equation*}
  \rho_0\ddot{u}(X, t) = \int_{\mc{H}(X)} \left[ \sN{T}[\xi][X, t] - \sN{T}[-\xi][X+\xi, t] \right] {\rm d}\xi ,
\end{equation*}
\begin{equation}
  - \rho_0 u_0 \omega^2(X,t) {\rm e}^{{\rm i}(kX-\omega t)} = \int_{\mc{H}(X)} \frac{\s{\omega}[\xi]}{\omega_0} \Big[ S_{(m)}(X) + S_{(m)}(X+\xi) \Big] \bigg(\frac{\sN{Y}^2}{\xi^2}\bigg)^{m-1} \frac{\sN{Y}}{\xi^2} {\rm d}\xi .
  \label{eqn:dispersion4}
\end{equation}
Due to the nature of the finite deformation model, this equation is clearly nonlinear and dependent on $X$, $t$, and $u_0$. Angular frequency of $X=0$ and $t=0$ is calculated using a Hookean classical model.
The dispersion curves are shown in \cref{fig:dispersion} for two different $u_0$. The periodic zero angular frequencies evident in all members of the material model are indicative of rank deficiency as discussed by \citet{belytschko2000unified}. As expected, for a very small displacement field $u_0$ the predicted behavior of all different members of the material class are the same.

\section{Discussion}
\label{sec:discussion}
In this section, the connection between the stability analysis presented in this work and a stability condition developed by \citet{silling2017stability} is discussed. Additionally, possible stabilization methods and the source of the instability are presented. 

\subsection{Silling's material stability test} 
\label{sec:sillingtest}
Using energy minimization, \citet{silling2017stability} derived a condition for stability
\begin{equation}
  {\rm d} \s{T} \bigcdot {\rm d} \s{Y} > 0 .
  \label{eqn:sillingTest}
\end{equation}
The equilibrium is stable if \cref{eqn:sillingTest} holds for a given displacement field $\s{Y}$. 

In the context of the finite deformation theory, using \cref{eqn:TC}, 
\begin{align*}
  {\rm d} \sN{T}_k = \s{\omega}[\bs{\xi}] \frac{\xi_p \, \xi_q}{|\bs{\xi}|^{(2m+2)}} L^{-1}_{pqij} \Big[ & {\rm d} S_{(m)ij} \, |\s{Y}|^{2m-2} Y_k + S_{(m)ij} 2(m-1) \, |\s{Y}|^{2m-4} Y_k Y_l \, {\rm d} Y_l + S_{(m)ij} \, |\s{Y}|^{2m-2} {\rm d} Y_k \Big],
\end{align*}
\begin{align}
  {\rm d} \s{T} \bigcdot {\rm d} \s{Y} = \int_{\mc{H}} & \Bigg\{ \s{\omega}[\bs{\xi}] \frac{\xi_p \, \xi_q}{|\bs{\xi}|^{(2m+2)}} L^{-1}_{pqij} \bigg[ {\rm d} S_{(m)ij} \, |\s{Y}|^{2m-2} Y_k \, {\rm d} Y_k + S_{(m)ij} 2(m-1) \, |\s{Y}|^{2m-4} Y_l \, {\rm d} Y_l Y_k {\rm d} Y_k + S_{(m)ij} \, |\s{Y}|^{2m-2} \, {\rm d} Y_k {\rm d} Y_k \bigg] \Bigg\} \, {\rm d} \bxi.
  \label{eqn:dTdotdY}
\end{align}
Specializing to the case of uniform deformations and perturbing a single material point, due to symmetry and uniform loading, the stress state does not change for the point itself (cf. \cref{rmk:0}), i.e. ${\rm d} S_{(m)ij} = 0$. Thus, using \cref{eqn:dTdotdY},
\begin{align*}
  {\rm d} \s{T} \bigcdot {\rm d} \s{Y} = \int_{\mc{H}} \s{\omega}[\bs{\xi}] \frac{\xi_p \, \xi_q}{|\bs{\xi}|^{(2m+2)}} L^{-1}_{pqij} \, S_{(m)ij} \bigg[ & 2(m-1) \, |\s{Y}|^{2m-4} Y_l \, {\rm d} Y_l Y_k \, {\rm d} Y_k + |\s{Y}|^{2m-2} \, {\rm d} Y_k \, {\rm d} Y_k \bigg] \, {\rm d} \bxi.
\end{align*}
Similar to \cref{sec:staticEquilibrium}, for a pure hydrostatic loading as in \cref{eqn:hydro}, using the constitutive relation in \cref{eqn:constitutive}, and limiting the perturbation to a single coordinate, e.g. ${\rm d} Y_1 = -\eps$ and ${\rm d} Y_2 = {\rm d} Y_3 = 0$, we obtain
\begin{equation}
  {\rm d} \s{T} \bigcdot {\rm d} \s{Y} = -\kappa \eps^2 |1+a|^{2m-2} a L^{-1}_{pqii} \int_{\mc{H}} \s{\omega}[\bs{\xi}] \frac{\xi_p \, \xi_q}{|\bs{\xi}|^{4}}  \bigg[ 2(m-1) \frac{\xi_l \, \xi_l}{|\bs{\xi}|^2} + 1 \bigg] \, {\rm d} \bxi.
  \label{eqn:altS}
\end{equation}
Note the similarity between \cref{eqn:altS} and \cref{eqn:Gamma}. The similarity is expected as the energy minimization approach and the Jacobian matrix method are rooted to the same idea. The explicit Jacobian matrix approach allows an investigation into the cause of the instability.

\subsection{Possibility of stabilization}
\label{sec:stabilize}
Several stabilization methods were suggested for the original correspondence model, whose instability is mostly associated with nonuniform deformations. \citet{silling2017stability} added a penalty term to the material model, similar to the method of \citet{littlewood2010simulation}, which resists deviation from homogeneous deformations. \citet{breitzman2018bond} and \citet{chen2018bond} proposed two different bond-level deformation gradients to reduce the influence of the average nodal integration and increase the impact of each bond deformation in evolving its own force state. \citet{chowdhury2019modified} somehow removed the instability in some problems by using an arbitrary split of a family of bonds into multiple sub-horizon regions. It comes from the mathematical formulation of the original model that the averaging integral operator allows the cancellation of nonuniform parts of deformation, which is the main reason behind the unstable behavior. All the mentioned stabilization techniques aim to eliminate the spurious unphysical zero-energy modes that appear in a particle discretization of the model. The finite deformation correspondence theory, in contrast, has a different issue as homogeneous deformations can cause instability in its continuum level. Therefore, none of the mentioned stabilization methods can address the observed unstable behavior, and any potentially stabilization method would look {\em artificial}. 

Another idea is to choose the $m$-member of the material class depending on the nature of the problem. For example, using $m>m_{\rm cr}$ when dealing with tension, and $m<m_{\rm cr}$ for compression. This idea, however, is doubted to be practical. First of all, as discussed in \cref{sec:numerical}, instabilities are not limited to pure hydrostatic loadings and can happen in a variety of loading scenarios. It will be tedious, if possible, to distinguish all the stable $m$-models for every possible loading condition. In addition, significant inconsistencies would emerge in dealing with history-dependent problems, such as plasticity.

Pursuing a stabilization method is not recommended for a number of reasons. First, the finite deformation correspondence theory was supposed to naturally improve the stability of the original model. However, a different type of material instability is inherent in the model, which appears even for uniform deformations. Second, the finite deformation theory is computationally more expensive than the original model. In addition, the recent theory is limited to isotropic materials. Therefore, seeking its stabilization seems pointless, and using an existing stabilized version of the original model is more favorable.

\subsection{Source of the instability}
\label{sec:source}
The origin of the problem is in fact rooted to the underlying classical model. It is well understood that the local Seth-Hill strain family suffers from the observed instability \citep{ball1976convexity}, i.e. $m<0$ group has tensile issues and $m>0$ shows compression difficulties. However, the problem is considerably milder in the classical model such that extremely larger tension or compression is required to cause an unstable behavior. The nonlocal integral operator in the correspondence theory largely exaggerates the problem, which results in small compressive or tensile deformations causing instability. The unstable behavior imposes significant difficulties in application of the finite deformation peridynamic model to any practical engineering problem. 

\begin{remark}
\label{rmk:sandia}
Recently, \citet{behzadinasab2019sandia}, within the context of Sandia Fracture Challenge 2017 \citep{kramer2019sandia}, used the $m=1$ model to predict the deformations and failure of an additively manufacture metal under tensile loading. While the application of the theory resulted in a qualitatively good blind simulation results, the load-carrying capacity of the structure was underestimated and an early failure was predicted. Although tension was the dominant loading mode of the problem, compression and shear states were inevitable in the material, due to a complicated geometry. While the superiority of the peridynamic Green-Lagrange strain model in dealing with tensile loading was showed (see \cref{sec:singular}), even instability in a small region can destroy the accuracy of the macroscopic response.
\end{remark}

\section{Conclusions}
\label{sec:conclusions}
The stability of a recently proposed peridynamic finite deformation correspondence model \citep{foster2018generalized} is studied in a detailed examination. The finite deformation theory was developed by using a family of nonlinear strain measure and reported to improve the stability of the original correspondence model \citep{silling2007peridynamic} whose instability is attributed to nonuniform deformations and particle discretizations. The offered analysis, however, shows that while the newer model has made some improvements over the original theory (cf. \cref{sec:singular}), it fails to maintain stability for homogeneous deformations in the mathematical foundation of the theory. The observed spurious unphysical deformation modes do not emerge because of particle discretization, which, actually interestingly, enhances the stability (cf. \cref{rmk:5}). The unstable behavior is originated from the underlying local theory, in a less severe manner, and exaggerated largely through the nonlocal integral operator. The identified instability can significantly lower the reliability of the model in practical problems where complicated loading scenarios occur (cf. \cref{rmk:sandia}). Since the recent model fails to naturally improve the stability of the original correspondence theory, is more computationally costly, and is limited to isotropic materials, its stabilization is not recommended. A stabilized version of the original model would be more desirable (see \cref{sec:stabilize}).

\section*{Acknowledgments}
\label{sec:acknowledge}
The authors are grateful for the financial support provided by the AFOSR MURI Center for Materials Failure Prediction through Peridynamics: Project NO. ONRBAA12-020. M. Behzadinasab also acknowledges the fellowship funding by The University of Texas at Austin.

\appendix

\renewcommand*{\thesection}{\appendixname~\Alph{section}}

\section{Simplified form of the shape tensor in the bulk of a Two-dimensional material} 
\label{sec:AppendixA}
Constraining the influence state to spherical, i.e. $\s{\omega}_s\an\bxi = \s{\omega}_s(|\bs{\xi}|)$, and noting symmetry, integration of all the odd combinations of indices in the shape tensor results in zero. The shape tensor is then evaluated to be \citep[Section 3.2]{foster2018generalized}
\begin{align*}
  L_{ijkl} = & \ \Big[ \delta_{ij}\delta_{kl}(1-\delta_{ik}) + \delta_{ik}\delta_{jl}(1-\delta_{ij}) + \delta_{il}\delta_{jk}(1-\delta_{ij}) \Big] \int_{\mc{H}} \mkern-4mu \s{\omega}_s\an{\bxi} \frac{\xi_1^2\xi_2^2}{|\bxi|^4} \dxi + \delta_{ij}\delta_{jk}\delta_{kl} \int_{\mc{H}} \s{\omega}_s\an{\bxi} \frac{\xi_1^4}{|\bxi|^4} \dxi .
\end{align*}
These integrals are then computed in 2D:
\begin{equation*}
  \int_{\mc{H}} \s{\omega}_s\an{\bxi} \frac{\xi_1^2\xi_2^2}{|\bxi|^4} \dxi = \frac{V_{\s{\omega}}}{8} , \qquad \int_{\mc{H}} \s{\omega}_s\an{\bxi} \frac{\xi_1^4}{|\bxi|^4} \dxi = \frac{3 V_{\s{\omega}}}{8},  
\end{equation*}
where $V_{\s{\omega}}$ is the weighted volume 
\begin{equation*}
  V_{\s{\omega}} = \int_{\mc{H}} \s{\omega}_s\an{\bxi} \dxi,
\end{equation*}
therefore,
\begin{equation*}
  L_{ijkl} = \frac{V_{\s{\omega}}}{8} (\delta_{ij}\delta_{kl} + \delta_{ik}\delta_{jl} + \delta_{il}\delta_{jk}).
\end{equation*}

Due to its isotropic nature, $\mbb{L}^{-1}$ can be written as:
\begin{equation}
  L_{klmn}^{-1} = A\delta_{kl}\delta_{mn} + B(\delta_{km}\delta_{ln} + \delta_{kn}\delta_{lm}), 
  \label{eqn:Linv2D}
\end{equation}
where $A$ and $B$ are unknown constants to be find through the following process:
\begin{equation*}
  L_{ijkl}L_{klmn}^{-1} = \frac{1}{2}(\delta_{im}\delta_{jn} + \delta_{in}\delta_{jm}) , \nonumber 
\end{equation*}
\begin{equation*}
  \frac{V_{\s{\omega}}}{8} (\delta_{ij}\delta_{kl} + \delta_{ik}\delta_{jl} + \delta_{il}\delta_{jk}) (A\delta_{kl}\delta_{mn} + B\delta_{km}\delta_{ln} + B\delta_{kn}\delta_{lm}) = \frac{1}{2}(\delta_{im}\delta_{jn} + \delta_{in}\delta_{jm}) , \nonumber 
\end{equation*}
\begin{equation*}
  A(2\delta_{ij}\delta_{mn} + \delta_{ij}\delta_{mn} + \delta_{ij}\delta_{mn}) + 2B(\delta_{ij}\delta_{mn} + \delta_{im}\delta_{jn} + \delta_{in}\delta_{jm}) = \frac{4}{V_{\s{\omega}}}(\delta_{im}\delta_{jn} + \delta_{in}\delta_{jm}) , \nonumber 
\end{equation*}
\begin{equation*}
  2(2A+B)\delta_{ij}\delta_{mn} + 2B(\delta_{im}\delta_{jn} + \delta_{in}\delta_{jm}) = \frac{4}{V_{\s{\omega}}}(\delta_{im}\delta_{jn} + \delta_{in}\delta_{jm}) .
\end{equation*}
Comparing the two sides of the equation 
\begin{equation*}
  2A+B = 0, \qquad B = \frac{2}{V_{\s{\omega}}},
\end{equation*}
hence, 
\begin{equation*}
  A = -\frac{1}{V_{\s{\omega}}}.
\end{equation*}
Substituting these back into \cref{eqn:Linv2D}, the inverted shape tensor is obtained 
\begin{equation*}
  L_{klmn}^{-1} = \frac{1}{V_{\s{\omega}}}(2\delta_{km}\delta_{ln} + 2\delta_{kn}\delta_{lm} - \delta_{kl}\delta_{mn}) .
\end{equation*}

\section{Simplifying \cref{eqn:2Ddfs}} 
\label{sec:appendixB}
In this section, \cref{eqn:2Ddfs} is explicitly simplified by breaking it into different parts.
{\small
\begin{align*}
  & \sum_{\mbf{X^J}\in\mc{H^I}} \frac{\xi^{\I\J}_k\xi^{\I\J}_l(4\xi^{\I\J}_i\xi^{\I\J}_j - \xi^{\I\J}_p\xi^{\I\J}_p\delta_{ij}) (4\xi^{\I\J}_r\xi^{\I\J}_s - \xi^{\I\J}_q\xi^{\I\J}_q\delta_{rs}) }{\left|\bs{\xi^{\I\J}}\right|^{8}} \Delta V^{\J} \nonumber \\
  & \qquad = 16 \sum_{\mbf{X^J}\in\mc{H^I}} \frac{\xi^{\I\J}_i\xi^{\I\J}_j\xi^{\I\J}_k\xi^{\I\J}_l\xi^{\I\J}_r\xi^{\I\J}_s}{\left|\bs{\xi^{\I\J}}\right|^{8}} - 4\delta_{rs} \sum_{\mbf{X^J}\in\mc{H^I}} \frac{\xi^{\I\J}_i\xi^{\I\J}_j\xi^{\I\J}_k\xi^{\I\J}_l}{\left|\bs{\xi^{\I\J}}\right|^{6}} - 4\delta_{ij} \sum_{\mbf{X^J}\in\mc{H^I}} \frac{\xi^{\I\J}_k\xi^{\I\J}_l\xi^{\I\J}_r\xi^{\I\J}_s}{\left|\bs{\xi^{\I\J}}\right|^{6}} + \delta_{ij} \delta_{rs} \sum_{\mbf{X^J}\in\mc{H^I}} \frac{\xi^{\I\J}_k\xi^{\I\J}_l}{\left|\bs{\xi^{\I\J}}\right|^{4}} , \nonumber \\
  & \qquad = 16 \bigg[ \gamma_{ijklrs} b \frac{\pi}{8} \ln N + \delta_{ki} \delta_{ij}\delta_{jr}\delta_{rs}\delta_{sl} b \Big(4 + \frac{5\pi}{8} \ln N\Big) \bigg] \notag \\
  & \qquad \quad - 4 \delta_{rs} \bigg[ \Big( \delta_{ij}\delta_{kl}(1-\delta_{ik}) + \delta_{ik}\delta_{jl}(1-\delta_{ij}) + \delta_{il}\delta_{jk}(1-\delta_{ij}) \Big) b\frac{\pi}{4}\ln N + \delta_{ki}\delta_{ij}\delta_{jl} b\Big(4 + \frac{3\pi}{4}\ln N\Big) \bigg] \nonumber \\
  & \qquad \quad - 4 \delta_{ij} \bigg[ \Big( \delta_{rs}\delta_{kl}(1-\delta_{rk}) + \delta_{rk}\delta_{sl}(1-\delta_{rs}) + \delta_{rl}\delta_{sk}(1-\delta_{rs}) \Big) b\frac{\pi}{4}\ln N + \delta_{kr}\delta_{rs}\delta_{sl} b\Big(4 + \frac{3\pi}{4}\ln N\Big) \bigg] + \delta_{ij}\delta_{rs}\delta_{kl} b(4 + \pi \ln N) , \\
  & \qquad = 4b \Bigg\{ \gamma_{ijklrs} \frac{\pi}{2} \ln N + \delta_{ki} \delta_{ij}\delta_{jr}\delta_{rs}\delta_{sl} \Big(16 + \frac{5\pi}{2} \ln N\Big) + \delta_{ij}\delta_{rs}\delta_{kl} (1 + \frac{\pi}{4} \ln N) \nonumber \\
  & \qquad \qquad \ \ - \delta_{rs} \bigg[ \Big( \delta_{ij}\delta_{kl}(1-\delta_{ik}) + \delta_{ik}\delta_{jl}(1-\delta_{ij}) + \delta_{il}\delta_{jk}(1-\delta_{ij}) \Big) \frac{\pi}{4}\ln N + \delta_{ki}\delta_{ij}\delta_{jl} \Big(4 + \frac{3\pi}{4}\ln N\Big) \bigg] \nonumber \\
  & \qquad \qquad \ \ - \delta_{ij} \bigg[ \Big( \delta_{rs}\delta_{kl}(1-\delta_{rk}) + \delta_{rk}\delta_{sl}(1-\delta_{rs}) + \delta_{rl}\delta_{sk}(1-\delta_{rs}) \Big) \frac{\pi}{4}\ln N + \delta_{kr}\delta_{rs}\delta_{sl} \Big(4 + \frac{3\pi}{4}\ln N\Big) \bigg] \Bigg\} .
\end{align*}
}
%
{\small 
\begin{align*}
  \sum_{\mbf{X^J}\in\mc{H^I}} \frac{(4\xi^{\I\J}_i\xi^{\I\J}_j - \xi^{\I\J}_n\xi^{\I\J}_n\delta_{ij})}{\left|\bs{\xi^{\I\J}}\right|^{4}} \bigg[ 2(m-1)\frac{\xi_k^{\I\J} \xi^{\I\J}_l}{\left|\bs{\xi^{\I\J}}\right|^2} + \delta_{kl} \bigg] \Delta V^{\J} 
  & = 8(m-1) \sum_{\mbf{X^J}\in\mc{H^I}} \frac{\xi^{\I\J}_i\xi^{\I\J}_j\xi^{\I\J}_k\xi^{\I\J}_l}{\left|\bs{\xi^{\I\J}}\right|^{6}} + 4\delta_{kl} \sum_{\mbf{X^J}\in\mc{H^I}} \frac{\xi^{\I\J}_i\xi^{\I\J}_j}{\left|\bs{\xi^{\I\J}}\right|^{4}} - 2(m-1)\delta_{ij} \sum_{\mbf{X^J}\in\mc{H^I}} \frac{\xi^{\I\J}_k\xi^{\I\J}_l}{\left|\bs{\xi^{\I\J}}\right|^{4}} - \delta_{ij} \delta_{kl} \sum_{\mbf{X^J}\in\mc{H^I}} \frac{1}{\left|\bs{\xi^{\I\J}}\right|^{2}} , \nonumber \\
  & = 8(m-1) \bigg[ \Big( \delta_{ij}\delta_{kl}(1-\delta_{ik}) + \delta_{ik}\delta_{jl}(1-\delta_{ij}) + \delta_{il}\delta_{jk}(1-\delta_{ij}) \Big) b\frac{\pi}{4}\ln N + \delta_{ki}\delta_{ij}\delta_{jl} b\Big(4 + \frac{3\pi}{4}\ln N\Big) \bigg] \nonumber \\
  & \qquad + 4\delta_{kl}\delta_{ij}b\left(4 + \pi \ln N\right) - 2(m-1)\delta_{ij}\delta_{kl}b\left(4 + \pi \ln N\right) - \delta_{ij}\delta_{kl}b\left(8 + 2\pi \ln N\right) , \nonumber \\
  & = 8(m-1)b \bigg[ \Big( \delta_{ij}\delta_{kl}(1-\delta_{ik}) + \delta_{ik}\delta_{jl}(1-\delta_{ij}) + \delta_{il}\delta_{jk}(1-\delta_{ij}) \Big) \frac{\pi}{4}\ln N + \delta_{ki}\delta_{ij}\delta_{jl} \Big(4 + \frac{3\pi}{4}\ln N\Big) \bigg] \nonumber \\
  & \qquad + 2(2-m)b \delta_{ij}\delta_{kl} \left(4 + \pi \ln N\right)
\end{align*}
}

\section{Simplifying \cref{eqn:2Ddfss}} 
\label{sec:appendixC}
Using \cref{eqn:gamma} (note that indices can be either 1 or 2, for a two-dimensional setup)
\begin{align*}
  \delta_{ij}\delta_{rs}\gamma_{ijklrs} = 3\delta_{kl} .
\end{align*}
Then, for the hydrostatic loading condition of \cref{eqn:hydro2Ds,eqn:hydro2De}, \cref{eqn:2Ddfss} is simplified 
\begin{align*}
  \p{f^\I_k}{y^\I_l} 
  = & - \frac{4}{\pi^2 N^2 (\Delta X)^2} (1+a)^{4m-2} \kappa \delta_{ij} \delta_{rs} \Bigg\{ \gamma_{ijklrs} \frac{\pi}{2} \ln N + \delta_{ki}\delta_{ij}\delta_{jr}\delta_{rs}\delta_{sl} \Big(16 + \frac{5\pi}{2} \ln N\Big) + \delta_{ij}\delta_{rs}\delta_{kl} (1 + \frac{\pi}{4} \ln N) \nonumber \\
  & \qquad \qquad \qquad \qquad \qquad \qquad \quad - \delta_{rs} \bigg[ \Big( \delta_{ij}\delta_{kl}(1-\delta_{ik}) + \delta_{ik} \delta_{jl}(1-\delta_{ij}) + \delta_{il}\delta_{jk}(1-\delta_{ij}) \Big) \frac{\pi}{4}\ln N + \delta_{ki}\delta_{ij}\delta_{jl} \Big(4 + \frac{3\pi}{4}\ln N\Big) \bigg] \nonumber \\
  & \qquad \qquad \qquad \qquad \qquad \qquad \quad - \delta_{ij} \bigg[ \Big( \delta_{rs}\delta_{kl}(1-\delta_{rk}) + \delta_{rk}\delta_{sl}(1-\delta_{rs}) + \delta_{rl}\delta_{sk}(1-\delta_{rs}) \Big) \frac{\pi}{4}\ln N + \delta_{kr}\delta_{rs}\delta_{sl} \Big(4 + \frac{3\pi}{4}\ln N\Big) \bigg] \Bigg\} \nonumber \\
  & - \frac{4}{\pi N (\Delta X)^2} (1+a)^{2m-2} \kappa \frac{(1+a)^{2m} - 1}{2m} \delta_{ij} \Bigg\{ 4(m-1)\bigg[ \Big(\delta_{ij}\delta_{kl}(1-\delta_{ik}) + \delta_{ik}\delta_{jl}(1-\delta_{ij}) + \delta_{il}\delta_{jk}(1-\delta_{ij})\Big) \frac{\pi}{4}\ln N \nonumber \\
  & \qquad \qquad \qquad \qquad \qquad \qquad \qquad \qquad \qquad \qquad \quad + \delta_{ki}\delta_{ij}\delta_{jl} \Big(4 + \frac{3\pi}{4}\ln N\Big) \bigg] + (2-m) \delta_{ij} \delta_{kl} \left(4 + \pi \ln N\right) \Bigg\} , \nonumber \\
  = & - \frac{4}{\pi^2 N^2 (\Delta X)^2} (1+a)^{4m-2} \kappa \Bigg\{ 3 \delta_{kl} \frac{\pi}{2} \ln N + \delta_{kl} \Big(16 + \frac{5\pi}{2} \ln N\Big) + 4\delta_{kl} (1 + \frac{\pi}{4} \ln N) - 4\delta_{kl} \Big( 4 + \pi \ln N \Big) \Bigg\} \nonumber \\
  & - \frac{4}{\pi N (\Delta X)^2} (1+a)^{2m-2} \kappa \frac{(1+a)^{2m} - 1}{m} \Bigg\{ 4(m-1)\bigg[\delta_{kl} \frac{\pi}{4}\ln N + \delta_{kl} \Big(4 + \frac{3\pi}{4}\ln N\Big)\bigg] + 2(2-m)\delta_{kl} \left(4 + \pi \ln N\right) \Bigg\} , \nonumber \\
  = & - \frac{4 (4 + \pi \ln N)}{\pi^2 N^2 (\Delta X)^2} \kappa (1+a)^{4m-2} \left[ 1 + 2 \pi N (1 - (1+a)^{-2m}) \right] \delta_{kl} .
\end{align*} 

\bibliography{refs}

\end{document}